\documentclass[12pt]{article}
\usepackage{fullpage,graphicx,psfrag,amsmath,amsfonts,verbatim}
\usepackage[small,bf]{caption}\usepackage{amssymb}
\usepackage{enumitem}
\usepackage{tikz}
\usepackage{listings}
\usepackage{authblk}
\usepackage{subfiles}
\usepackage{url}
\usepackage{textcomp}
\usepackage{booktabs}
\usepackage{algorithm}
\usepackage{algpseudocode}
\usepackage{caption}
\usepackage{subfig}
\usepackage{hyperref}

\usepackage{xcolor}
\hypersetup{
    colorlinks=true,     
    linkcolor=blue,      
    citecolor=red,       
    urlcolor=magenta     
}

\newtheorem{theorem}{Theorem}

\lstdefinestyle{mystyle}{
      backgroundcolor=\color{backcolour},
      keywordstyle=\color{magenta},
      numberstyle=\tiny\color{codegray},
      stringstyle=\color{codepurple},
      basicstyle=\ttfamily\small,
      breakatwhitespace=false,
      breaklines=true,
      captionpos=t,
      keepspaces=true,
      numbers=left,
      numbersep=5pt,
      showspaces=false,
      showstringspaces=false,
      showtabs=false,
      tabsize=2,
}

\newcommand{\ones}{\mathbf 1}
\newcommand{\reals}{{\mbox{\bf R}}}

\newcommand{\symm}{{\mbox{\bf S}}}  

\newcommand{\Tr}{\mathop{\bf Tr}}
\newcommand{\diag}{\mathop{\bf diag}}

\newcommand{\prox}{\mathbf{prox}}

\newcommand{\cl}{{\mathop {\bf cl}}}

\newcommand{\dom}{\mathop{\bf dom}} 

\newcommand{\sign}{\mathop{\bf sign}}

\newcommand{\ie}{{\it i.e.,}}

\newcommand{\BEAS}{\begin{eqnarray*}}
\newcommand{\EEAS}{\end{eqnarray*}}
\newcommand{\BEA}{\begin{eqnarray}}
\newcommand{\EEA}{\end{eqnarray}}
\newcommand{\BEQ}{\begin{equation}}
\newcommand{\EEQ}{\end{equation}}
\newcommand{\BIT}{\begin{itemize}}
\newcommand{\EIT}{\end{itemize}}
\newcommand{\BNUM}{\begin{enumerate}}
\newcommand{\ENUM}{\end{enumerate}}

\newcommand{\BEQN}{\begin{equation*}}
\newcommand{\EEQN}{\end{equation*}}

\newcounter{algorithmctr}
\renewcommand{\thealgorithmctr}{\arabic{algorithmctr}}
   {\mbox{}\\*[\parskip]\begin{minipage}{\linewidth}%
       \refstepcounter{algorithmctr}\begin{list}{}{%
       \setlength{\rightmargin}{0\linewidth}%
       \setlength{\leftmargin}{.05\linewidth}}%
       \rmfamily\small
       \item[]{\setlength{\parskip}{0ex}\hrulefill\par%
        \nopagebreak{\bfseries\textsf{Algorithm \thealgorithmctr~}}}}%
   {{\setlength{\parskip}{-1ex}\nopagebreak\par\hrulefill\\*[2ex]\par}%
   \end{list}\end{minipage}}

\title{Projections onto Spectral Matrix Cones}
\author[]{Daniel Cederberg \qquad \qquad Stephen Boyd \\ Stanford University}
\date{November 2, 2025}

\begin{document}
\maketitle

\begin{abstract}
Semidefinite programming is a fundamental problem class in convex optimization, 
but despite recent advances in solvers, solving large-scale semidefinite
programs remains challenging. Generally the matrix functions involved
are \emph{spectral} or \emph{unitarily invariant}, \ie~they depend only 
on the eigenvalues or singular values of the matrix. This paper investigates 
how \emph{spectral matrix cones} --- cones defined from epigraphs and perspectives
of spectral or unitarily invariant functions --- can be used to enhance 
first-order conic solvers for semidefinite programs. Our main result
shows that projecting a matrix can be reduced to projecting its eigenvalues 
or singular values, which we demonstrate can be done at a negligible cost
compared to the eigenvalue or singular value decomposition itself.
We have integrated support for spectral matrix cone projections into the
Splitting Conic Solver (SCS). Numerical experiments show that SCS with this
enhancement can achieve speedups of up to an order of magnitude for solving
semidefinite programs arising in experimental design, robust principal component 
analysis, and graph partitioning.
\end{abstract}
      
\newpage
\tableofcontents
\newpage
\section{Introduction}

\paragraph{Conic convex optimization.} Many general-purpose solvers for convex
optimization use a standard form similar to
\BEQ \label{e:standard-form}
\begin{array}{ll}
\mbox{minimize} & c^T x \\
\mbox{subject to} & Ax + s = b \\
& s \in K, \\
\end{array}
\EEQ 
where the decision variables are $x \in \reals^n$ and $s \in \reals^{m}$,
and the problem data are the matrix $A \in \reals^{m \times n}$, the vectors $b
\in \reals^m$ and $c \in \reals^n$, and the nonempty closed convex cone $K
\subseteq \reals^m$. While manually reformulating convex problems into this standard
form can be error-prone and tedious, modeling languages such as YALMIP
\cite{Lofberg2004} and CVXPY \cite{Diamond2016} facilitate the process by taking
a high-level problem description and transforming it to the standard form
through a process called \emph{canonicalization}.

After canonicalization, the problem is expressed using the so-called
\emph{standard cones}, which are the nonnegative cone, the second-order cone,
the positive semidefinite cone, and the three-dimensional exponential and power
cones \cite{boyd2004, ODonoghue2016}. Canonicalizing a problem into a form based
on the standard cones often requires the introduction of numerous auxiliary
variables and constraints \cite{aps2025mosek}, ultimately resulting in a larger
problem formulation. For certain problem classes, such as linear programs (LPs)
and quadratic programs (QPs), the larger formulation does not typically
degrade solver performance significantly, as the added structure is usually
sparse and many LP and QP solvers are good at exploiting sparsity. However, for
other problem classes, such as semidefinite programming, the increase in problem
size can degrade the solver performance. In many cases, smaller and simpler
conic formulations could be obtained if solvers supported a broader set of cones
beyond the standard ones. This is the motivation behind the recent
\emph{interior-point solver} Hypatia \cite{Coey2022, Coey2022a, Coey2023}, which
supports a much broader class of cones than standard solvers.

Interior-point solvers like Mosek \cite{mosekSolver} or Clarabel
\cite{Goulart2024} are frequently the preferred choice for many optimization
problems, but scaling up the problem size with those solvers can sometimes be
challenging. This limitation is particularly pronounced for semidefinite
programs (SDPs), where it is not uncommon for problems involving a dense matrix
variable of dimension 150 to be beyond the computational reach of interior-point
methods on a standard laptop. An alternative approach that scales better to
large-scale problems is to employ first-order solvers such as ADMM-based SCS
\cite{ODonoghue2016, Odonoghue21} and COSMO \cite{Garstka2021}, or PDHG-based
methods \cite{applegate2021, lin2025}. In every iteration, those methods project
onto the cone $K$, where the projections for all standard cones are either
available in closed form or can be computed efficiently through an iterative
process \cite{Parikh2014}. Furthermore, recent research on more specialized
cones has focused on designing efficient algorithms for computing projections
with respect to a generalized distance measure, replacing the traditional
\emph{Euclidean} projection \cite{Chao2018, Jiang2022, Cederberg2025}. By
designing efficient algorithms for projecting onto cones beyond the standard
cones, first-order conic solvers can be applied to simpler and smaller conic
formulations of problems compared to if only standard cones were used,
potentially resulting in improved solver performance. 

\paragraph{Spectral matrix cones.} In this paper, we consider a class of cones
that we will refer to as \emph{spectral matrix cones}. These cones are based on
the notion of \emph{spectral functions} \cite{Davis1957, Friedland1981,
Lewis1995, Lewis1996}, which are functions $F:\symm^n \to \reals \cup
\{\infty\}$  satisfying 
\[
F(X) = F(U X U^T) 
\] 
for all orthogonal matrices $U \in \reals^{n \times n}$. In other words, the
function $F$ is a spectral function if the value $F(X)$ only depends on the
unsorted eigenvalues of $X$. A simple example is $F(X) = \Tr(X) = \sum_{i=1}^n
\lambda_i(X)$.

Given a convex spectral function $F$, we define a set $K_F \subseteq \reals
\times \reals_{++} \times \symm^n$ to be a \emph{spectral matrix cone} if it can
be expressed as the closure of the epigraph of the perspective of $F$:
\[
K_F \triangleq \cl \{(t, v, X) \in \reals \times
\reals_{++} \times \symm^n \: | \: v F(X/v) \leq t\}.    
\] 
The set $K_F$ is a closed convex cone, since the perspective function 
$(X, v) \mapsto v F(X/v)$ of $F$ is convex and positively homogeneous for $v >
0$ \cite{boyd2004}. This construction of a convex cone from the epigraph of a
perspective function is well known \cite[\S 5]{Rockafellar1970} and underpins
the common folklore that any convex optimization problem can be expressed as a
conic linear program \cite{Nemirovski2006}.

A first example of a spectral matrix cone is obtained from the log-determinant
function $F(X) = - \log \det X$ with domain $\dom F = \symm^{n}_{++}$. This
spectral function induces a so-called \emph{log-determinant cone}, which can be
used to obtain much simpler and smaller conic formulations of problems involving
the log-determinant function compared to its standard canonicalization based on
the positive semidefinite cone. (We discuss this further in \S
\ref{sec:potential advantage}.)
  
\paragraph{Spectral vector cones and projections onto spectral matrix cones.}
The advantage of the spectral matrix cones is that we will be able to reduce the
problem of projecting a \emph{matrix} to the problem of projecting the
\emph{eigenvalues} of the matrix. This simplification is based on a classic
result \cite[Corollary 1]{Davis1957} on spectral functions saying that a convex
function $F:\symm^n \to \reals \cup \{\infty\}$ is a spectral function if and
only if there exists a convex function $f:\reals^n \to \reals \cup \{\infty\}$
that is \emph{symmetric} (meaning that $f(x) = f(Px)$ for all permutation
matrices $P \in \reals^{n \times n}$) such that
\[
    F(X) = f(\boldsymbol{\lambda}(X)),    
\] 
where $\boldsymbol{\lambda}(X) = (\lambda_1(X), \lambda_2(X), \dots
\lambda_n(X)) \in \reals^{n}$  is the vector of eigenvalues of $X$ in
nonincreasing order. We use this characterization to prove that projecting
$(\bar{t}, \bar{v}, \bar{X}) \in \reals \times \reals \times \symm^{n}$ onto
$K_F$ can be reduced to projecting $(\bar{t}, \bar{v},
\boldsymbol{\lambda}(\bar{X})) \in \reals \times \reals \times \reals^n$ onto
the cone
\[
K_f \triangleq \cl \{(t, v, x) \in \reals \times
\reals_{++} \times \reals^n \: | \: v f(x/v) \leq t\},   
\] 
where $f$ is the symmetric convex function corresponding to $F$. To distinguish
$K_f$ from $K_F$ we will refer to it as a \emph{spectral vector cone}. In other
words, to project onto $K_F$ we first compute an eigenvalue decomposition and
then project onto $K_f$. We demonstrate that for several spectral matrix cones,
projecting onto $K_f$ can be done using a structure-exploiting implementation of
Newton's method with an iteration cost of order $\mathcal{O}(n)$, or by sorting
an $n$-dimensional vector at a cost of order $\mathcal{O}(n \log n)$. Since the cost of
the eigenvalue decomposition is of order $\mathcal{O}(n^3)$, we expect it to
dominate the cost for the projection onto $K_F$. (In practice, as we demonstrate
in \S\ref{sec:ablation}, the time to compute the eigenvalue decomposition is
often more than two or three orders of magnitude greater than the time to project onto the
spectral vector cone $K_f$.)

\paragraph{Consequences of efficient projections.}
To demonstrate the importance of our results, we have integrated the new
projection routines for spectral matrix cones into the first-order conic solver
SCS \cite{ODonoghue2016, Odonoghue21}. Numerical experiments demonstrate that
this enhanced version of SCS, with support for projections onto spectral matrix
cones, is often an order of magnitude faster than the default version for
certain types of SDPs. The performance gains can be attributed to two
factors: (1) faster iterations due to smaller eigenvalue decompositions (we
explain this further in \S\ref{sec:potential advantage}), and (2) a potential
reduction in the number of iterations required for convergence.

\paragraph{Extension to unitarily invariant functions on $\reals^{m \times n}$.}
While our main focus in this paper is cones defined in terms of spectral
functions $F:\symm^n \to \reals \cup \{\infty \}$ whose value $F(X)$ only
depends on the \emph{eigenvalues} of $X \in \symm^n$, our results trivially
extend to cones defined in terms of functions $F:\reals^{m \times n} \to \reals
\cup \{\infty\}$ whose value $F(X)$ only depends on the \emph{singular values}
of $X \in \reals^{m \times n}$. In this setting, we consider functions
$F:\reals^{m \times n} \to \reals \cup \{\infty\}$ that are \emph{unitarily
invariant}, meaning that 
\[
  F(X) = F(UXV)  
\] 
for all unitary matrices $U \in \reals^{m \times m}$ and $V \in \reals^{n \times n}$.
A well known result \cite{Lewis1995} is that a convex function $F:\reals^{m
\times n} \to \reals \cup \{\infty\}$ where $m \geq n$ is unitarily invariant if
and only if there exists a convex function $f:\reals^n \to \reals \cup
\{\infty\}$ that is \emph{absolutely symmetric} (meaning that $f(x) = f(|Px|)$
for all permutation matrices $P \in \reals^{n \times n}$, where $|\cdot|$ is
elementwise absolute value)  such that 
\[
F(X) = f(\boldsymbol{\sigma}(X)),  
\] 
where $\boldsymbol{\sigma}(X) = (\sigma_1(X), \dots, \sigma_n(X)) \in
\reals^{n}$ is the vector of singular values of $X$ in nonincreasing order. This
characterization of unitarily invariant functions can be used to derive results
that are very similar to the corresponding results for cones defined in terms of
spectral functions. For brevity we will sometimes omit the details. 

\paragraph{Additional related work.} 
Optimization algorithms that rely on projections, along with the development of
efficient algorithms for computing projections onto structured sets, have been
central themes in the optimization literature for decades, as reflected by a
steady string of papers over the years \cite{gubin1967method, bregman1967,
dykstra1983, bauschke1996projection, duchi2008efficient, fadili2010total,
nemeth2010project, condat2016fast, bauschke2018projecting, perez2020filtered,
briceno2024projection, roth2025n, Luxenberg25}. In particular, the projections
onto the standard cones are available in closed form \cite[\S 6]{Parikh2014},
with the exceptions of the exponential and power cones. The current
state-of-the-art for projecting onto the exponential cone relies on an iterative
algorithm and was recently described in \cite{Friberg2023}. An iterative
algorithm for projecting onto the power cone is described in
\cite{hien2015differential}.

The idea of incorporating customized routines within first-order
methods for projecting onto convex sets beyond the standard cones is briefly
mentioned in \cite{Garstka2021}. However, no specific examples of cones 
that might benefit from such tailored projection algorithms are provided.

Our main results are conceptually related to several works that employ a
\emph{general transfer principle}, which states   that various properties of 
function or sets in $\symm^n$ can be transferred to corresponding properties 
of functions or sets in $\reals^n$ \cite{Lewis1995, Lewis1996, Parikh2014, Daniilidis2014}.
In particular, in \cite{Daniilidis2008} and \cite[Appendix A]{Lewis2008}, they study
\emph{theoretical} properties of projections onto so-called \emph{spectral
sets}, which are sets of the form $\{X \in \symm^n \: | \:
\boldsymbol{\lambda}(X) \in C\}$ where $C \subseteq \reals^n$ is a symmetric
set, meaning that for every $x \in C$ and any permutation matrix $P \in
\reals^{n \times n}$ it holds that $Px \in C$. If $K_F$ is a spectral matrix
cone with associated spectral vector cone $K_f$ (according to the definition in
this paper), it can equivalently be expressed as 
\[
K_F = \{(t, v, X) \in \reals \times \reals \times \symm^n \: | \: 
(t, v, \boldsymbol{\lambda}(X)) \in K_f \}.
\] 
This formulation reveals similarities to the structure of spectral sets.
However, while the existing literature on spectral sets primarily focuses on
theoretical properties, the potential for spectral matrix
cones to enhance conic first-order solvers for semidefinite programming remains
unexplored.

\paragraph{Outline.} The rest of the paper is organized as follows. In \S
\ref{sec:potential advantage}, we present two examples that demonstrate the
potential advantages of incorporating spectral matrix cones within a first-order
conic solver. In \S \ref{sec:spectral matrix cone projections} we give our main
results, along with several pairs of spectral vector and matrix cones. We then
discuss how to efficiently project onto the spectral vector cones in \S
\ref{sec:spectral vector cone projections}. Numerical results for SCS using the
spectral matrix cone projections are given in \S \ref{sec:numerical
experiments}, followed by conclusions in \S \ref{sec:conclusions}.

\section{Motivating examples}
\label{sec:potential advantage}
In this section we give two examples of problems where spectral matrix cones
offer smaller and more compact conic formulations compared to the standard
canonicalization based on the standard cones. The first cone is induced by the
log-determinant function $F(X) = - \log \det X$, $\dom F = \symm^{n}_{++}$,
resulting in the \emph{log-determinant cone} 
\BEQ \label{e:logdet-cone}
K_{\text{logdet}} \triangleq \cl \big( 
\{(t, v, X) \in \reals \times \reals_{++} \times \symm^n_{++}  \: | \:
- v \log \det(X/v) \leq t \} \big).
\EEQ 
The second cone is induced by the nuclear norm $F(X) = \|X \|_*$, 
$\dom F = \reals^{m \times n}$, resulting in the \emph{nuclear norm cone}
\BEQ \label{e:nuc-cone}
K_{\text{nuc}} \triangleq \{(t, X) \in \reals \times \reals^{m \times n}  \: | \:
 \| X \|_* \leq t \} .   
\EEQ 
Since the nuclear norm is positively homogeneous, its perspective function $(v,
X) \mapsto v F(X/v)$ is independent of the perspective variable $v$. We
therefore consider the nuclear norm cone as a subset of $\reals \times \reals^{m
\times n}$ instead of $\reals \times \reals \times \reals^{m \times n}$. (For a
few other cones we will also use this convention that positively homogeneous
functions induce cones through their epigraphs without applying the
perspective function.)

\subsection{Log-determinant cone}
The log-determinant function has applications in many fields, including
computational geometry, statistics, and control \cite{Vandenberghe1998}. To
demonstrate how problems involving this function are canonicalized, we consider
the problem 
\BEQ \label{e:logdetcone-motivation 1}
\begin{array}{ll}
\mbox{minimize} & \Tr(SX) - \log \det X
\end{array}
\EEQ
with variable $X \in \symm^n$ and problem data $S \in \symm^{n}_{++}$. This is
the maximum likelihood estimation problem of the covariance matrix of a Gaussian
random vector, where $X$ is the inverse covariance matrix and $S$ is the
sample covariance \cite{boyd2004}. While this problem has the analytical
solution $X = S^{-1}$, we use it to demonstrate how the log-determinant cone can
give a more compact conic formulation compared to if only standard cones were
used for the canonicalization. The same reasoning applies to other problems
without analytical solution involving the log-determinant function.

The standard canonicalization of \eqref{e:logdetcone-motivation 1} is based on the 
fact that for a given $X \succ 0$, the value of $\log \det X$ is equal to the optimal
value of \cite[\S 6]{aps2025mosek}
\[
\begin{array}{ll}
\mbox{maximize} & \sum_{i=1}^n \log Z_{ii}\\
\mbox{subject to} & \begin{bmatrix} X & Z \\ Z^T & \diag(Z) \end{bmatrix} \succeq 0 \\
& Z \text{ lower triangular}, \\
\end{array}  
\] 
with variable $Z \in \reals^{n \times n}$. This allows us to 
rewrite \eqref{e:logdetcone-motivation 1} as 
\BEQ \label{e:logdetcone-motivation 2}
\begin{array}{ll}
\mbox{minimize} & \Tr(S X) - \sum_{i=1}^n \log Z_{ii} \\
\mbox{subject to} & \begin{bmatrix} X & Z \\ Z^T & \diag(Z)
\end{bmatrix} \succeq 0 \\
& Z \text{ lower triangular},
\end{array}
\EEQ 
with variables $X \in \symm^{n}$ and $Z \in \reals^{n \times n}$.
Formulation \eqref{e:logdetcone-motivation 2} can then be expressed in the
standard form \eqref{e:standard-form} using $n$ exponential cones (one for each
log-term in the objective) and a positive semidefinite cone of matrices with dimension $2n$.

To canonicalize \eqref{e:logdetcone-motivation 1} using the log-determinant
cone, we first write it as 
\[
\begin{array}{ll}
\mbox{minimize} & \Tr(SX) + t \\
\mbox{subject to} &  -\log \det X \leq t,
\end{array}
\]
with variables $t \in \reals$ and $X \in \symm^n$, which is equivalent 
to 
\BEQ \label{e:logdetcone-motivation 3}
\begin{array}{ll}
\mbox{minimize} & \Tr(SX) + t \\ 
\mbox{subject to} & v = 1 \\
&  (t, v, X) \in K_{\text{logdet}},
\end{array}
\EEQ 
with variables $t \in \reals, \: v \in \reals$, and $X \in \symm^n$.
In other words, this canonicalization requires only one log-determinant cone 
with a matrix variable of dimension $n$.

First-order methods applied to \eqref{e:logdetcone-motivation 2} must project
onto the positive semidefinite cone of matrices with dimension $2n$ in every
iteration. In contrast, when applied to \eqref{e:logdetcone-motivation 3},
first-order methods must project onto $K_{\text{logdet}}$, which we will show
is dominated by the cost of computing the eigenvalue decomposition of a
symmetric matrix of dimension $n$. Since the cost of computing an eigenvalue
decomposition scales cubically with the dimension of the matrix, the projection
onto the log-determinant cone can potentially be $8$ times faster. A first-order
method applied to \eqref{e:logdetcone-motivation 3} can therefore have
significantly faster iterations than a first-order method applied to
\eqref{e:logdetcone-motivation 2}.

\subsection{Nuclear norm cone}
The nuclear norm $\|X \|_* = \sum_{i=1}^n \sigma_i(X)$, where $X \in \reals^{m
\times n}$ and $\sigma_1(X) \geq \sigma_{2}(X) \geq \dots \geq \sigma_n(X)$ are
the singular values of $X$, commonly appears in convex heuristics for rank
minimization \cite{Fazel2002}. (Throughout this paper we will assume that $m
\geq n$.) To demonstrate how this function is canonicalized, we consider the
problem 
\BEQ \label{e:nuclearcone-motivation 1}
\begin{array}{ll}
\mbox{minimize} & \| X \|_* \\
\mbox{subject to} & \| S \|_1 \leq \mu \\
& X + S = M,
\end{array}
\EEQ 
with variables $X \in \reals^{m \times n}$ and $S \in \reals^{m \times n}$,
and problem data $M \in \reals^{m \times n}$ and $\mu > 0$. This problem is 
known as \emph{robust principal component analysis} and has been 
proposed for recovering a low rank matrix from measurements $M$ that have been
corrupted by sparse noise $S$ \cite{Candes2011}.

The standard canonicalization of \eqref{e:nuclearcone-motivation 1} is based on
the fact that for a given $X \in \reals^{m \times n}$, the value of $\| X \|_*$
is equal to the optimal value of \cite[\S 6]{aps2025mosek}
\[
\begin{array}{ll}
\mbox{minimize} & (1/2)(\Tr(U) + \Tr(V))\\
\mbox{subject to} & \begin{bmatrix} U & X^T \\ X & V \end{bmatrix} \succeq 0,
\end{array}  
\] 
with variables $U \in \reals^{n \times n}$ and $V \in \reals^{m \times m}$. 
This allows us to rewrite \eqref{e:nuclearcone-motivation 1}
as 
\BEQ \label{e:nuclearcone-motivation 2}
\begin{array}{ll}
\mbox{minimize} & (1/2)(\Tr(U) + \Tr(V))\\
\mbox{subject to} & \begin{bmatrix} U & X^T \\ X & V \end{bmatrix} \succeq 0 \\
& \| S\|_1 \leq \mu \\
& X + S = M, 
\end{array}  
\EEQ
with variables  $U \in \reals^{n \times n}, \: V \in \reals^{m \times m}, \: 
X \in \reals^{m \times n}$, and $S \in \reals^{m \times n}$. Formulation 
\eqref{e:nuclearcone-motivation 2} can then be expressed in the standard form
\eqref{e:standard-form} using $(1 + 2mn)$ nonnegative cones and a positive
semidefinite cone of matrices with dimension $m + n$. 

To canonicalize \eqref{e:nuclearcone-motivation 1} using the nuclear norm cone,
we recognize that it can be written as 
\BEQ \label{e:nuclearcone-motivation 3}
\begin{array}{ll}
\mbox{minimize} & t \\
\mbox{subject to} & \| S\|_1 \leq \mu \\
& X + S = M \\
& (t, X) \in K_{\text{nuc}},
\end{array}  
\EEQ
with variables $t \in \reals$, $X \in \reals^{m \times n}$, and $S \in
\reals^{m \times n}$. 

First-order methods applied to \eqref{e:nuclearcone-motivation 2} must project
onto the positive semidefinite cone of matrices with dimension $m + n$ in every iteration, resulting
in a cost of order $\mathcal{O}((m+n)^3)$. In contrast, when applied to
\eqref{e:nuclearcone-motivation 3}, first-order methods must project onto
$K_{\text{nuc}}$, which we will show is dominated by the cost of 
computing the reduced singular value decomposition (SVD) of an $m \times n$
matrix. Since the reduced SVD of an $m \times n$ matrix with $m \geq n$ can be
computed at a cost of order $\mathcal{O}(m n^2)$, a first-order method applied
to \eqref{e:nuclearcone-motivation 3} can have significantly faster
iterations than a first-order method applied to \eqref{e:nuclearcone-motivation
2}.

\section{Spectral matrix cone projections}
\label{sec:spectral matrix cone projections}
In this section we present the main results that allow us to project onto a
spectral matrix cone by first computing either the eigenvalue decomposition or
the singular value decomposition of the matrix, followed by projecting onto the
corresponding spectral vector cone. We also present several pairs of spectral
vector and matrix cones.

\subsection{Main results}
In the following result we denote the projection of $(\bar{t}, \bar{v}, \bar{X})
\in \reals \times \reals \times \symm^n$ onto the spectral matrix cone $K_F$ as
$\Pi_{K_F}(\bar{t}, \bar{v}, \bar{X}) \in \reals \times \reals \times \symm^n$,
and the projection of $(\bar{t}, \bar{v}, \bar{\lambda }) \in \reals \times
\reals \times \reals^n$ onto the spectral vector cone $K_f$ as
$\Pi_{K_f}(\bar{t}, \bar{v}, \bar{\lambda}) \in \reals \times \reals \times
\reals^n$.
\begin{theorem}
\label{e:main-thm}
Let $f:\reals^n \to \reals \cup \{\infty\}$ be a symmetric convex function
corresponding to the spectral function $F:\symm^n \to \reals \cup \{\infty\}$.
Consider the projection of $(\bar{t}, \bar{v}, \bar{X}) \in \reals \times \reals
\times \symm^n$ onto $K_F$, where $\bar{X}$ has eigenvalue decomposition
$\bar{X} = U \diag(\bar{\lambda}) U^T$ with $\bar{\lambda} \in \reals^n$. Then 
\[
\begin{split}
\Pi_{K_F}(\bar{t}, \bar{v}, \bar{X})_1 & = 
\Pi_{K_f}(\bar{t}, \bar{v}, \bar{\lambda})_1 \\
\Pi_{K_F}(\bar{t}, \bar{v}, \bar{X})_2 & = 
\Pi_{K_f}(\bar{t}, \bar{v}, \bar{\lambda})_2 \\   
\Pi_{K_F}(\bar{t}, \bar{v}, \bar{X})_3 & = 
U\diag(\Pi_{K_f} (\bar{t}, \bar{v}, \bar{\lambda})_3) U^T,
\end{split} 
\] 
where $\Pi_{K_F}(\bar{t}, \bar{v}, \bar{X})_1 \in \reals, \: \Pi_{K_F}(\bar{t},
\bar{v}, \bar{X})_2 \in \reals, \: \Pi_{K_F}(\bar{t}, \bar{v}, \bar{X})_3 \in
\symm^n$ and $\Pi_{K_f}(\bar{t}, \bar{v}, \bar{\lambda})_1 \in \reals, \:
\Pi_{K_f}(\bar{t}, \bar{v}, \bar{\lambda})_2 \in \reals, \: \Pi_{K_f}(\bar{t},
\bar{v}, \bar{\lambda})_3 \in \reals^n$ are the first, second and last
components of the projections $\Pi_{K_F}(\bar{t}, \bar{v}, \bar{X})$ and
$\Pi_{K_f}(\bar{t}, \bar{v}, \bar{\lambda})$, respectively. 
\end{theorem}
\emph{Proof.} The projection $\Pi_{K_F}(\bar{t}, \bar{v}, \bar{X})$ is
the solution of 
\BEQ \label{e:proof 1}
\begin{array}{ll}
\mbox{minimize} & (1/2) (t - \bar{t})^2 + (1/2) (v - \bar{v})^2 +
(1/2) \| X - \bar{X} \|_F^2 \\
\mbox{subject to} & (t, v, X) \in K_F,
\end{array} 
\EEQ 
with variable $(t, v, X) \in \reals \times \reals \times \symm^n$. It is
easy to verify that $(t, v, X) \in K_F$ if and only if $(t, v,
\boldsymbol{\lambda}(X)) \in K_f$. Furthermore, $\| X - \bar{X} \|_F^2 \geq  \|
\boldsymbol{\lambda}(X) - \boldsymbol{\lambda}(\bar{X}) \|_2^2$ for all $X \in
\symm^n$ \cite[\S 7.4]{Horn2013}. It follows that a lower bound on the optimal
value of \eqref{e:proof 1} is given by optimal value of 
\BEQ \label{e:proof 2}
\begin{array}{ll}
\mbox{minimize} & (1/2) (t - \bar{t})^2 + (1/2) (v - \bar{v})^2 + 
(1/2) \| \lambda - \bar{\lambda} \|_2^2 \\
\mbox{subject to} & (t, v, \lambda) \in K_f,
\end{array} 
\EEQ
with variable $(t, v, \lambda) \in \reals \times \reals \times \reals^n$.
The optimal solution of \eqref{e:proof 2} is $(t^\star, v^\star,
\lambda^\star) \triangleq \Pi_{K_f}(\bar{t}, \bar{v},
\bar{\lambda})$. The point $(t^\star, v^\star, U \diag(\lambda^\star)U^T)$ is
feasible in \eqref{e:proof 1} and attains the same objective value as the
optimal value of \eqref{e:proof 2}, which is a lower bound of the optimal value
of \eqref{e:proof 1}, and this point must therefore be optimal in \eqref{e:proof
1}. $\blacksquare$

We can state a similar result for a spectral matrix cone defined in terms of a 
unitarily invariant function.
\begin{theorem} \label{e:main-thm-singular}
Let $f:\reals^n \to \reals \cup \{\infty \}$ be an absolutely symmetric convex
function corresponding to the unitarily invariant function $F:\reals^{m\times n}
\to \reals \cup \{\infty\}$ where $m \geq n$. Consider the projection of
$(\bar{t}, \bar{v}, \bar{X}) \in \reals \times \reals \times \reals^{m \times
n}$ onto $K_F$, where $\bar{X}$ has reduced singular value decomposition $\bar{X} = U
\diag(\bar{\sigma}) V^T$ with $U \in \reals^{m \times n}, \: V \in \reals^{n \times n}$, and
singular values $\bar{\sigma} \in \reals^n$. Then 
\[
\begin{split}
\Pi_{K_F}(\bar{t}, \bar{v}, \bar{X})_1 & =
 \Pi_{K_f}(\bar{t}, \bar{v}, \bar{\sigma})_1 \\
\Pi_{K_F}(\bar{t}, \bar{v}, \bar{X})_2 & = 
\Pi_{K_f}(\bar{t}, \bar{v}, \bar{\sigma})_2 \\   
\Pi_{K_F}(\bar{t}, \bar{v}, \bar{X})_3 & = 
U \diag(\Pi_{K_f}(\bar{t}, \bar{v}, \bar{\sigma})_3) V^T.
\end{split}
\]
\end{theorem}
The proof of Theorem \ref{e:main-thm-singular} is very similar to the
proof of Theorem \ref{e:main-thm} and is therefore omitted.

\subsection{Examples of spectral cone pairs}
\label{sec:Examples}
By choosing different symmetric or absolutely symmetric convex functions, we can
derive several pairs of spectral vector and matrix cones.

\paragraph{Log-determinant cone.} Consider the symmetric convex function 
$f(x) = - \sum_{i=1}^n \log x_i$, $\dom f = \reals^{n}_{++}$. The corresponding
spectral function is $F(X) = -\log \det X$, $ \dom F = \symm^{n}_{++}$.
The associated spectral vector cone is the \emph{logarithmic cone}, given by
\[
K_{\text{log}}  \triangleq \cl  
\{(t, v, x) \in \reals \times \reals_{++} \times \reals^n_{++}  \: | \:
-  \sum_{i=1}^n v \log(x_i/v) \leq t \}.    
\] 
In Appendix \ref{sec:appendix-topology-log} we show that it can be expressed as 
\[
K_{\text{log}} =   \{(t, v, x) \in \reals \times \reals_{++} \times 
\reals^n_{++}  \: | \:
-  \sum_{i=1}^n v \log(x_i/v) \leq t \}
\cup \: (\reals_+ \times \{0\} \times \reals_+^n).
\]
This expression will later simplify the derivation of the projection onto
$K_{\text{log}}$. The associated spectral matrix cone is the
log-determinant cone, defined in \eqref{e:logdet-cone}.

\paragraph{Nuclear norm cone.} Consider the absolutely symmetric convex function 
$f(x) = \|x \|_1$, $\dom f = \reals^n$. The corresponding spectral function is 
the nuclear norm $F(X) = \|X \|_*$, $\dom F = \reals^{m \times n}$ where we
assume $m \geq n$. The associated spectral vector cone is the \emph{$\ell_1$-norm
cone}, given by 
\[
K_{\ell_1} \triangleq \{(t, x) \in \reals \times \reals^n \: | \: \|x \|_1 \leq t\}.  
\]
The associated spectral matrix cone is the nuclear norm cone, defined in 
\eqref{e:nuc-cone}.

\paragraph{Trace-inverse cone.} Consider the symmetric convex function 
$f(x) = \sum_{i=1}^n 1/x_i$, $\dom f = \reals^{n}_{++}$. The corresponding spectral
function is $F(X) = \Tr(X^{-1})$, $\dom F = \symm^{n}_{++}$. The associated 
spectral vector cone is the \emph{inverse cone}, given by
\[
K_{\text{inv}} \triangleq \cl \{(t, v, x) \in \reals \times \reals_{++}
\times \reals^n_{++} \: | \: v^2 \sum_{i=1}^n 1/x_i \leq t \}.    
\]
In Appendix \ref{sec:appendix-topology-inverse} we show that it can be expressed as 
\[
K_{\text{inv}} =   \{(t, v, x) \in \reals \times \reals_{++}
\times \reals^n_{++} \: | \: v^2 \sum_{i=1}^n 1/x_i \leq t \}
\cup \: (\reals_+ \times \{0\} \times \reals_+^{n}).
\]
This expression will later simplify the derivation of the projection onto 
$K_{\text{inv}}$.
The associated spectral matrix cone is the \emph{trace-inverse cone}, given by 
\[
K_{\text{TrInv}} \triangleq \cl \{(t, v, X) \in \reals \times \reals_{++}
\times \symm^n_{++} \: | \: v^2 \Tr(X^{-1}) \leq t \}.    
\]

\paragraph{Entropy cone.} Consider the symmetric convex function $f(x) =
\sum_{i=1}^n x_i \log x_i$, $\dom f = \reals^{n}_{+}$ with the convention that
$0 \log 0 = 0$. The corresponding spectral function is the (negative)
von-Neumann entropy $F(X) = \sum_{i=1}^n \lambda_i(X) \log \lambda_i(X)$, $\dom
F = \symm^{n}_+$. (For $X \in \symm^{n}_{++}$ it holds that $F(X) = \Tr(X \log
X)$.) The associated spectral vector cone is the \emph{vector entropy cone},
given by 
\[
K_{\text{vEnt}} \triangleq \cl \{ (t, v, x) \in \reals \times \reals_{++}
\times \reals^{n}_{+} \: | 
\: \sum_{i=1}^n x_i \log (x_i/v) \leq t\}.
\]
This cone is also known as the \emph{relative entropy cone} 
\cite{Chandrasekaran2017}.
In Appendix \ref{sec:appendix-topology-entropy} we show that it can be expressed as 
\[
K_{\text{vEnt}} =  \{ (t, v, x) \in \reals \times \reals_{++} \times 
\reals^{n}_{+} \: | \: \sum_{i=1}^n x_i \log (x_i/v) \leq t\}
\cup \: (\reals_+ \times \reals_+ \times \{0\}^n).
\]
This expression will later simplify the derivation of the projection onto
$K_{\text{vEnt}}$. The associated spectral matrix cone is the \emph{matrix
entropy cone}, given by
\[
K_{\text{mEnt}} \triangleq \cl \{ (t, v, X) \in \reals \times \reals_{++} 
\times \symm^{n}_{++} \: 
| \: \Tr(X \log (X/v)) \leq t \}.    
\]

\paragraph{Root-determinant cone.} Consider the symmetric convex function $f(x)
= -\prod_{i=1}^n x_i^{1/n}$, $\dom f = \reals^{n}_+$. The corresponding spectral
function is the root-determinant function $F(X) = -(\det(X))^{1/n}$, $\dom F =
\symm^{n}_{+}.$ The associated spectral vector cone is the \emph{geometric mean
cone}, given by
\[
K_{\text{geomean}}  \triangleq \{(t, x) \in \reals \times \reals^n_{+} 
\: | \: - \prod_{i=1}^n x_i^{1/n} \leq t\}.
\] 
The associated spectral matrix cone is the \emph{root-determinant cone}, given by
\[
K_{\det} \triangleq \{(t, X) \in \reals \times \symm^n_{+} \: | \: -( \det(X))^{1/n} \leq t \}.     
\]

\paragraph{Sum-of-largest-eigenvalues cone.}
For $x \in \reals^n$ we denote by $x_{[i]}$ the \emph{i}th largest component of
$x$, \ie~$x_{[1]} \geq x_{[2]} \geq \dots \geq x_{[n]}$. For any $k \in \{1, 2,
\dots, n\}$, consider the symmetric convex function $f(x) = \sum_{i=1}^k
x_{[i]}$, $\dom f = \reals^n$. The corresponding spectral function is the
sum-of-\emph{k}-largest eigenvalue function $F(X) = \sum_{i=1}^k
\boldsymbol{\lambda}_i(X)$, $\dom F = \symm^n.$ The associated spectral vector
cone is the \emph{sum-of-largest-entries cone}, given by
\[
K_{\text{vSum}} \triangleq \{ (t, x) \in \reals \times \reals^n \: | \: 
\sum_{i=1}^k x_{[i]} \leq t\}.    
\]
The associated spectral matrix cone is the \emph{sum-of-largest-eigenvalues
cone}, given by
\[
K_{\text{mSum}} \triangleq \{ (t, X) \in \reals \times \symm^n \: | \: 
\sum_{i=1}^k \boldsymbol{\lambda}_i(X) \leq t\}.    
\]

\paragraph{Dual cones.} The Moreau decomposition \cite[\S 2.5]{Parikh2014}
allows us to project onto the dual cone $K^*_F$ by projecting onto the closed
convex cone $K_F$. It holds that
\BEQ \label{e:Moreau} 
\Pi_{K_F}(t, v, X) - \Pi_{K^*_F}(-(t, v, X)) = (t, v, X),  
\EEQ 
where $\Pi_{K_F}(t, v, X)$ is the projection of $(t, v, X)$ onto $K_F$,
and $\Pi_{K^*_F}(-(t, v, X))$ is the projection of $-(t, v, X)$ onto $K^*_F$.
The dual cone of the epigraph of the perspective of a function is equal to the
epigraph of the perspective of the conjugate function, but with the epigraph and
perspective components swapped \cite[\S 14]{Rockafellar1970}:
\[
\begin{split}
K^*_f & =  \cl \{(t, v, x) \in \reals_{++} \times
\reals \times \reals^n \: | \: v \geq t f^*(-x/t)\} \\
K^*_F & =  \cl \{(t, v, X) \in \reals_{++} \times
\reals \times \symm^n \: | \: v \geq t F^*(-X/t)\}.
\end{split}
\] 
(Here $F^*(X) = \sup_{Y} \{ \Tr(XY) - F(Y) \}$ and $f^*(x) = \sup_y \{ x^T y -
f(y )\}$ are the Fenchel conjugates of $F$ and $f$, respectively.) By using a
classic result on spectral functions \cite[Theorem 2.3]{Lewis1996} saying that
$F^*(X) = f^*(\boldsymbol{\lambda}(X))$, we derive several dual cones which we
list in Appendix \ref{sec:appendix dual cones}. Deriving explicit expressions for
the dual cones is important, as they will later be used to verify optimality in
the spectral vector cone projection problems.

\section{Spectral vector cone projections}
\label{sec:spectral vector cone projections}
In the previous section we showed that we can reduce the problem of projecting
onto a spectral matrix cone to the problem of projecting onto the corresponding
spectral vector cone. In this section we discuss how to efficiently project onto
the spectral vector cones. First we present techniques for projecting onto the
$\ell_1$-norm cone and the sum-of-largest-entries cone, which require ad-hoc
analysis. We then present projections for the remaining cones using a more
systematic approach based on Newton's method. 

\subsection{Ad-hoc projections}
\label{sec:ad-hoc projection}
\paragraph{$\ell_1$-norm cone.} The projection of a point $(\bar{t}, \bar{x})
\in \reals \times \reals^n$ onto $K_{\ell_1}$ can be found by sorting the
entries of $\bar{x}$ in descending order of magnitude, resulting in a total cost
of order $\mathcal{O}(n \log n)$. The following algorithm can be derived using
the Moreau decomposition \eqref{e:Moreau} together with the fact that it is
known how to project onto the dual cone $K_{\ell_1}^*$ efficiently (see, for
example, \cite{Ding2014}). Let $\pi \in \reals^{n}$ be a permutation of $\{1,
\dots, n\}$ such that $|\bar{x}_{\pi(1)}| \geq |\bar{x}_{\pi(2)}| \geq \dots
\geq |\bar{x}_{\pi(n)}|.$ There exists a unique $k \in \{0, 1, 2, \dots, n\}$ such
that 
\[
|\bar{x}_{\pi(k)} | > \max \left\{ \frac{1}{k+1}\left(-\bar{t} + \sum_{i=1}^k 
|\bar{x}_{\pi(k)}|\right), 0 \right\} \geq |\bar{x}_{\pi(k+1)}|,     
\]
where we use the convention that $\bar{x}_{\pi(0)} = \infty$ and 
$\bar{x}_{\pi(n+1)} = \infty$. The projection $(t, x)$ of $(\bar{t}, \bar{x})$
onto $K_{\ell_1}$ is given by 
\[
\begin{split}
\hspace{4cm} t & =  \max \left\{ \frac{1}{k+1}\left(k\bar{t} + \sum_{i=1}^k 
|\bar{x}_{\pi(k)}|\right), \bar{t} \right\} \\
\hspace{4cm} x_i & = \max(|\bar{x}_i| - (t - \bar{t}), 0) \sign{\bar{x}_i},
\qquad i = 1, \dots, n.
\end{split}
\]
Note that the $x$-update is simply equal to $\prox_{(t - \bar{t}) \| \cdot \|_1}(\bar{x})$, 
\ie~the proximal operator of the (scaled) $\ell_1$-norm evaluated at $\bar{x}$.

When the $\ell_1$-norm cone projection is used for projecting onto a nuclear
norm cone, the vector $\bar{x}$ is the singular values of a matrix, and the
entries are therefore guaranteed to be sorted and nonnegative. This makes the
projection onto the $\ell_1$-norm cone extremely cheap. (In our experiments in
\S\ref{sec:ablation}, the cost of computing the SVD is several
orders of magnitude larger than the overhead associated with the projection onto
the $\ell_1$-norm cone.)

\paragraph{Sum-of-largest-entries cone.} Recently in \cite{Luxenberg25}, an
efficient algorithm for projecting onto sublevel sets of the
sum-of-$k$-largest function $f(x) = \sum_{i=1}^k x_{[i]}$ was derived. In
Appendix \ref{sec:appendix sum-of-largest} we extend their algorithm to project
a point $(\bar{t}, \bar{x}) \in \reals \times \reals^n$ onto the
sum-of-largest-entries cone $K_{\text{vSum}}$. The extended algorithm can be
summarized as follows. 
\BNUM 
\item Compute a permutation $\pi \in \reals^n$ of $\{1, \dots, n\}$ such that
  $\bar{x}_{\pi(1)} \geq \bar{x}_{\pi(2)} \geq \dots \geq \bar{x}_{\pi(n)}$ by
  sorting $\bar{x}$. Denote the sorted vector by $\bar{x}_s$.
\item Compute scalars $\eta$, $n_u$, $n_t$, $a_t$ by calling Algorithm
\ref{alg:sum-largest} (introduced below) with input $(\bar{t}, \bar{x}_s)$.
\item The projection of $(\bar{t}, \bar{x}_{s})$ onto $K_{\text{vSum}}$ is $(t_s, x_s)$
where $t_s = \bar{t} + \eta$ and 
\[
x_{s} = \big(\underbrace{(\bar{x}_{{s}})_1 - \eta, \dots, 
(\bar{x}_{{s}})_{n_u} - \eta}_{n_u}, \underbrace{a_t, \dots, a_t}_{n_t}, 
(\bar{x}_{{s}})_{n_u + n_t + 1}, \dots, (\bar{x}_{{s}})_{n}  \big).
\]
\item The projection of $(\bar{t}, \bar{x})$ onto $K_{\text{vSum}}$ is $(t, x)$ where $t = t_s$ and 
$x$ is obtained from $x_s$ by undoing the sort in step 1 using the permutation
$\pi$.
\ENUM 
The complexity of the algorithm is determined by the sort in the first step,
resulting in a complexity of order $\mathcal{O}(n \log n)$. This is the same
complexity as simply evaluating if $(\bar{t}, \bar{x})$ belongs to
$K_{\text{vSum}}$. We present more details on the algorithm in Appendix 
\ref{sec:appendix sum-of-largest}.

We should also point out that when the projection onto the
sum-of-largest-entries cone is used for projecting onto a
sum-of-largest-eigenvalues cone, the vector $\bar{x}$ is the eigenvalues of a
matrix, and it is therefore guaranteed to be sorted. In this case the first and
last steps above are not needed. 

\begin{algorithm}
\caption{Computing scalars for projecting $(\bar{t}, \bar{x}) \in \reals \times \reals^n$ onto
$K_{\text{vSum}}$, assuming $0 < k < n$ and that $\bar{x}$ is sorted.}\label{alg:sum-largest}
\begin{footnotesize}
\begin{algorithmic}[1]
\State \textbf{Input:} $\bar{t} \in \reals, \: \bar{x} \in \reals^n, \: k
\in \{1, 2, \dots, n-1\}$
\State \textbf{Initialize:} $n_u \gets k$, $\eta \gets 0$, $S \gets
\sum_{i=1}^{k} \bar{x}_{i}$, $a_u \gets \bar{x}_{n_u}$, $a_t \gets
\bar{x}_{n_u + 1}$, $t \gets \bar{t}$, $n_t \gets 0$
\While{$S > t$}
\State \textbf{if} $n_u = k$ \textbf{then} $r \gets 1$ \textbf{else}
$r \gets n_t /(k - n_u)$
\State \textbf{if} $n_u = k$ \textbf{then} $s_1 \gets a_u - a_t$
\textbf{else if} $n_u = 0$ \textbf{then} $s_1 \gets \infty$
\textbf{else} $s_1 \gets (a_u - a_t) /(r - 1)$
\State \textbf{if} $n_u + n_t = n$ or $n_t = 0$ \textbf{then} $s_2
\gets \infty$
\textbf{else} $s_2 \gets a_t - \bar{x}_{n_u + n_t + 1}$
\State $s_3 \gets (S - t) / (r(n_u + 1) + k - n_u)$
\State $s \gets \min(s_1, s_2, s_3)$
\State $\eta \gets \eta + sr$, $S \gets S - s(r n_u + k - n_u)$, $t \gets t_0 + \eta$
\State \textbf{if} $n_t > 0$ \textbf{then} $a_t \gets a_t - s$
\State \textbf{if} $s = s_1$ \textbf{then} $n_u \gets n_u - 1$
\State \textbf{if} $n_u > 0$ \textbf{then} $a_u \gets \bar{x}_{n_u} - \eta$
\State \textbf{if} $n_t = 0$ \textbf{then} $n_t \gets 2$ \textbf{else} $n_t \gets n_t + 1$ 
\EndWhile
\State $n_t \gets \max(n_t - 1, 0)$
\State \textbf{return} $\eta$, $n_u$, $n_t$, $a_t$
\end{algorithmic}
\end{footnotesize}
\end{algorithm}

\subsection{Systematic projections}

\subsubsection{Main idea}
For the remaining spectral vector cones in \S \ref{sec:Examples}, we compute the 
projection of $(\bar{t}, \bar{v}, \bar{x})$ onto $K_f$ by applying an iterative
method to solve 
\BEQ \label{e:proj-problem}
\begin{array}{ll}
\mbox{minimize} &  (1 / 2) ( t - \bar{t} )^2 + (1 / 2) ( v - \bar{v} )^2
 + (1 / 2) \| x - \bar{x} \|_2^2  \\
\mbox{subject to} & (t, v, x) \in K_f,
\end{array}
\EEQ 
with variables $t \in \reals$, $v \in \reals$, and $x \in \reals^n$.
The optimality conditions of \eqref{e:proj-problem}  are
\BEQ \label{e:proj-opt-cond}
\begin{split}
(t, v, x) & \in K_f, \qquad (\bar{t}, \bar{v}, \bar{x}) = 
(t, v, x) - (\lambda_t, \lambda_v, \lambda_x) \\
(\lambda_t, \lambda_v, \lambda_x) & \in K^*_f, \qquad t \lambda_t + 
v \lambda_v + \lambda_x^T x = 0,
\end{split}
\EEQ 
where $\lambda_t \in \reals$, $\lambda_v \in \reals$, and $\lambda_x \in \reals^n$ are 
Lagrange multipliers. These optimality conditions can be used to show that 
the negative dual cone is projected onto the origin, \ie~if $(\bar{t}, \bar{v},
\bar{x}) \in -K_f^*$, then its projection onto $K_f$ is 0. 

A spectral vector cone $K_f$ can be expressed as 
\[
K_f = \{(t, v, x) \in \reals \times \reals_{++} \times \reals^n \: | \: s_f(v, x) \leq t \}
\cup K_f^0 
\] 
for some convex function $s_f:\reals_{++} \times \reals^n \to \reals \cup \{\infty\}$ and set 
$K_f^0 \subseteq \reals^{n+2}$. (For example, for the logarithmic cone we have
$s_f(v, x) = -\sum_{i=1}^n v \log(x_i/v), \: \dom s_f = \reals^{n+1}_{++}$ and $K_f^0 = \reals_+ \times \{0\}
\times \reals^n_+$.) To solve \eqref{e:proj-problem} we first check if either
$(\bar{t}, \bar{v}, \bar{x}) \in K_f$ or $(\bar{t}, \bar{v}, \bar{x}) \in
-K^*_f$, or if the solution belongs to $K_f^0$ by using the optimality
conditions \eqref{e:proj-opt-cond}. If none of these cases apply, the projection
(\ie~the solution of \eqref{e:proj-problem}) can be found by solving 
\BEQ \label{e:proj-problem 2}
\begin{array}{ll}
\mbox{minimize} &  (1 / 2) ( t - \bar{t} )^2 +
(1 / 2) ( v - \bar{v} )^2 + (1 / 2) \| x - \bar{x} \|_2^2  \\
\mbox{subject to} & s_f(v, x) \leq t.
\end{array}
\EEQ 
To solve \eqref{e:proj-problem 2} we use the following two-step procedure.

\paragraph{First step.} At optimality of \eqref{e:proj-problem 2}, the 
inequality $s_f(v, x) \leq t$ is satisfied as an equality whenever $(\bar{t},
\bar{v}, \bar{x}) \not \in K_f$. We therefore do the substitution $t = s_f(v, x)$
and solve 
\BEQ \label{e:non-convex prob}
\begin{array}{ll}
\mbox{minimize} & h(v, x) = (1 / 2) ( s_f(v, x) - \bar{t} )^2 +
(1 / 2) ( v - \bar{v} )^2 + (1 / 2) \| x - \bar{x} \|_2^2  \\
\end{array}
\EEQ 
using Newton's method \cite[\S 9.5]{boyd2004}. A good starting point for
Newton's method is often available from the solution of the projection problem
in the previous iteration. Warmstarting Newton's method from this point often
makes the method extremely efficient, with only a handful of iterations
necessary for convergence (typically between two and five). However, there are
two subtleties with this approach:
\BNUM 
\item Problem \eqref{e:non-convex prob} is an unconstrained \emph{non-convex}
problem, and it occasionally happens that Newton's method converges to a
stationary point of the objective function that does not correspond to the
solution of the convex problem \eqref{e:proj-problem 2}. (Nevertheless, it is
possible to show that if Newton's method 
converges to a point $(v^\star, x^\star)$ with $s_f(v^\star, x^\star) > \bar{t}$,
then $(s_f(v^\star, x^\star), v^\star, x^\star)$ is the solution of
\eqref{e:proj-problem 2}.)
\item The domain of the function $s(v, x)$ includes the implicit constraints $v >
0$ and $x \succ 0$. This implicit constraint is not taken into account in Newton's method
(except in the line search), and we have observed
that this causes Newton's method to sometimes converge to the origin, even
though the origin is not the solution of \eqref{e:proj-problem 2}.
\ENUM
While both of these cases occasionally occur, we should point out that a
warmstarted Newton's method applied to \eqref{e:non-convex prob} finds the
solution of \eqref{e:proj-problem 2} in a vast majority of cases.  

\paragraph{Second step.} After applying Newton's method to \eqref{e:non-convex
prob} and obtaining a candidate solution $(v^\star, x^\star)$, we evaluate
whether $(s_f(v^\star, x^\star), v^\star, x^\star)$ solves \eqref{e:proj-problem
2} using the optimality conditions \eqref{e:proj-opt-cond}. If these conditions
are not met, we proceed by applying a primal-dual interior-point method (IPM)
directly to \eqref{e:proj-problem 2}. We have observed that a basic primal-dual
IPM, as the one described in \cite[\S 11.7]{boyd2004}, does not perform well,
and we have therefore implemented an enhanced version with an adaptive choice of
the centering parameter based on the progress made by the affine-scaling
direction \cite{Mehrotra1992, Vandenberghe2010}. We also include a higher-order
correction term \cite{Mehrotra1992}, a non-monotonic line search \cite[\S
15.6]{NoceWrig06}, and iterative refinement \cite{Higham2002}. We also noted
that scaling the point to be projected such that all its entries have magnitude
at most one, improved the performance. After the scaled point has been
projected, its projection is unscaled to obtain the projection of the original
point. (This is related to positive homogeneity of the Euclidean projection
operator onto a closed convex cone; see, for example, \cite{Ingram1991}.)

\subsubsection{Exploiting structure in Newton's method}
When applying Newton's method to \eqref{e:non-convex prob}, the search direction
$(\Delta v, \Delta x)$ is computed by solving the linear system 
\[
  \nabla^2 h(v, x) \begin{bmatrix} \Delta v \\ \Delta x \end{bmatrix} = - \nabla h(v, x).
\]
For each spectral vector cone in \S \ref{sec:Examples}, the Hessian $\nabla^2
h(v, x)$ has a structure that allows us to solve this linear system at a cost of
order $\mathcal{O}(n)$. For example, for the logarithmic cone, the Hessian is
$\nabla^2 h(v, x) = D(v, x) + w(v, x) w(v, x)^T$ with 
\[
\begin{split}
D(v, x)  = I + a  \diag(-a/v^2& + n/v - 2c/v,  v z_1^2, v z_2^2, \dots, v z_n^2)  \\
w(v, x)  & = \begin{bmatrix} - (a/v  + c) \\ v z \end{bmatrix}, 
\end{split}
\] 
where $z = 1 / x$ elementwise, $a = -\bar{t} - v \sum_{i=1}^n \log(x_i/v)$ and
$c = n - \sum_{i=1}^n \log(x_i/v)$. Linear systems with this
diagonal-plus-rank-one structure can be solved at a cost of order
$\mathcal{O}(n)$ \cite[Appendix C]{boyd2004}. 

\subsubsection{Exploiting structure in the interior-point method}
To solve \eqref{e:proj-problem 2} with an interior-point method, we add an 
epigraph variable $r \in \reals$ and target the formulation
\[
\begin{array}{ll}
\mbox{minimize} & r \\
\mbox{subject to} & f_i(t, v, x, r) \leq 0, \qquad i = 0, 1, 2 \\
\end{array}
\] 
where 
$f_0(t, v, x, r) = (1/2) ( t - \bar{t} )^2 + (1/2) ( v - \bar{v} )^2 +
(1/2) \| x - \bar{x} \|_2^2 - r$,
$f_1(t, v, x, r) = s(v, x) - t$ and $f_2(t, v, x, r) = - v$. Let $u = (t, v, x,
r)$ denote the primal variable and $z \in \reals^3$ the Lagrange multipliers for
the constraints $f_i(u) \leq 0, i = 0, 1, 2.$ The main cost of each iteration of 
the interior-point method is to solve linear systems of the form 
\BEQ \label{eq:IPM linear sys 1}
\begin{bmatrix}
H & B^T \\
B & - D 
\end{bmatrix}
\begin{bmatrix}
\Delta u \\
\Delta z
\end{bmatrix} = 
\begin{bmatrix}
b_1 \\ b_2
\end{bmatrix},
\EEQ 
for some right-hand side specified by $b_1 \in \reals^{n+3}$ and $b_2 \in
\reals^3$, and where $D \in \reals^{3 \times 3}$ is a diagonal matrix with
positive entries, $B \in \reals^{3 \times (n+3)}$ is the Jacobian of $f(u) =
(f_0(u), f_1(u), f_2(u))$, and $H = z_0 \nabla^2 f_0(u) + z_1 \nabla^2 f_1(u) +
z_2 \nabla^2 f_2(u)$ is the Hessian of the Lagrangian \cite[\S 8]{Wright1997}.
To solve \eqref{eq:IPM linear sys 1} we eliminate $\Delta z$ using the second
equation, resulting in the equation 
\[
(H + B^T D^{-1} B) \Delta u = b_1 + B^T D^{-1} b_2.
\] 
For each spectral vector cone in \S \ref{sec:Examples}, the Hessian $H$ of the
Lagrangian has a structure that allows us to solve this linear system at a cost
of order $\mathcal{O}(n)$. The efficient computation of the solution is based on
interpreting the system as a low-rank perturbation of a system with coefficient
matrix $H$, followed by applying the matrix inversion formula \cite[Appendix
C]{boyd2004}. However, it is well known that solving low-rank perturbed linear
systems using the matrix inversion formula can provide inaccurate solutions
\cite{Yip1986}. In our preliminary experiments, this inaccuracy caused the
interior-point method to struggle with finding a highly accurate solution of
\eqref{e:proj-problem 2}. To resolve this issue it was necessary to include
iterative refinement \cite{Higham2002}. 

\section{Numerical experiments}
\label{sec:numerical experiments}
In this section we examine the impact of incorporating spectral matrix cone
projections into the first-order conic optimization solver SCS \cite{ODonoghue2016,
Odonoghue21}. The new projections have been integrated into SCS
which is available at \url{https://github.com/cvxgrp/scs}. Code to run the
experiments below are available at 
\url{https://github.com/dance858/SpectralSCSExperiments}.
All the experiments have been carried out on a machine with a
13th Generation Intel® Core™ i7-1355U CPU running Ubuntu version 24.04.

To evaluate the performance of SCS with and without the spectral matrix cones,
we compare the following metrics in addition to the total solve time (we refer
to SCS enhanced with projections onto spectral matrix cones as
\texttt{SpectralSCS}, and the standard SCS as \texttt{SCS}):
\BNUM 
\item \emph{Total number of iterations}. To gain further insight into the
difference in total solve time, we show the total number of iterations required
for convergence.
\item \emph{Time per iteration}. In each iteration, both \texttt{SCS} and
\texttt{SpectralSCS} project onto a matrix cone (\ie~\texttt{SCS} projects onto
the positive semidefinite cone and \texttt{SpectralSCS} projects onto a spectral
matrix cone), and solve a linear system with cached factorization.
\texttt{SpectralSCS} can potentially have much faster cone projections by, for
example, computing a smaller eigenvalue decomposition compared to when the
standard canonicalization based on the positive semidefinite cone is used (an
example where this occurred was given in \S \ref{sec:potential advantage}). To show how
faster cone projections contribute to faster iterations, we show the average
time per iteration.
\item \emph{Time per matrix cone projection}. For \texttt{SCS} we show the time
required to project onto the positive semidefinite cone. For
\texttt{SpectralSCS}, we show the total time for the projection onto the
spectral matrix cone. (Later in \S \ref{sec:ablation}, we compare the time to
compute the eigenvalue decomposition with the time to project onto the spectral
vector cone for \texttt{SpectralSCS}.)
\ENUM
For each problem type and dimensions, we average the result over 5 problem
instances. Both \texttt{SCS} and \texttt{SpectralSCS} are used with the default
parameter settings, unless otherwise stated.

The implementation of the operator splitting algorithm that both \texttt{SCS} and 
\texttt{SpectralSCS} are based on depends on a dual scaling parameter 
$\sigma > 0$. The choice of $\sigma$ can have a large impact on the performance,
and by default this parameter is estimated adaptively using 
a heuristic that aims to balance the convergence rate of the primal and dual
residuals. Roughly speaking, this heuristic increases $\sigma$ if the primal
residual is much larger than the dual, and decreases it if the opposite is true.
However, for \texttt{SpectralSCS} and the log-determinant cone (but
not the other cones, or for \texttt{SCS}), the value of $\sigma$ sometimes seemed to
decrease even if the primal residual was much larger and was increasing. 
We therefore disabled the heuristic and used the default value $\sigma = 0.1$
for \texttt{SpectralSCS} on the experiments for the log-determinant cone.

\subsection{Log-determinant cone}

\subsubsection{Experimental design}
We consider the problem of computing the minimum-volume ellipsoid, centered at
the origin, that contains a given set of points $v_1, \dots, v_p \in \reals^n$.
This problem has applications in computational geometry and experimental design
\cite{boyd2004, Todd2016}. The ellipsoid is of the form $\{x \in \reals^n \: |
\: x^T W^\star x \leq 1\}$, where $W^\star \in \symm^n$ is the solution of
\BEQ \label{e:exp design}
\begin{array}{ll}
\mbox{minimize} & -\log \det W \\
\mbox{subject to} &  v_i^T W v_i \leq 1, \: i = 1, \dots, p,
\end{array}
\EEQ 
with variable $W \in \symm^n$ and problem data $v_i \in
\reals^n, \: i = 1, \dots, p$. To solve \eqref{e:exp design} with
\texttt{SpectralSCS} we formulate it as 
\[
\begin{array}{ll}
\mbox{minimize} & t \\
\mbox{subject to} &  v_i^T W v_i \leq 1, \: i = 1, \dots, p \\
& v = 1 \\
& (t, v, W) \in K_{\text{logdet}},
\end{array}
\]
with variables $t \in \reals, \: v \in \reals$, and $W \in \symm^n$.
The main cost of the matrix cone projection for \texttt{SCS} and
\texttt{SpectralSCS} is to compute the eigenvalue decomposition of a symmetric
matrix of dimension $2n$ and $n$, respectively.
(We should point out that the Cholesky factorization $W = L L^T$ allows a change
of variables, with lower-triangular $L \in \reals^{n \times n}$ as the new variable, 
to canonicalize \eqref{e:exp design} using only second-order and exponential cones.
However, we don't apply this canonicalization trick, since our goal is to compare 
\texttt{SCS} and \texttt{SpectralSCS}.) 

To generate the problem data we follow \cite{Coey2022}, choosing $p = 2n$ for different 
values of $n$, and generating $v_1, \dots, v_p$ with independent zero mean Gaussian
entries.

\begin{figure}[!htb]
\centering
\includegraphics[width=1\textwidth]{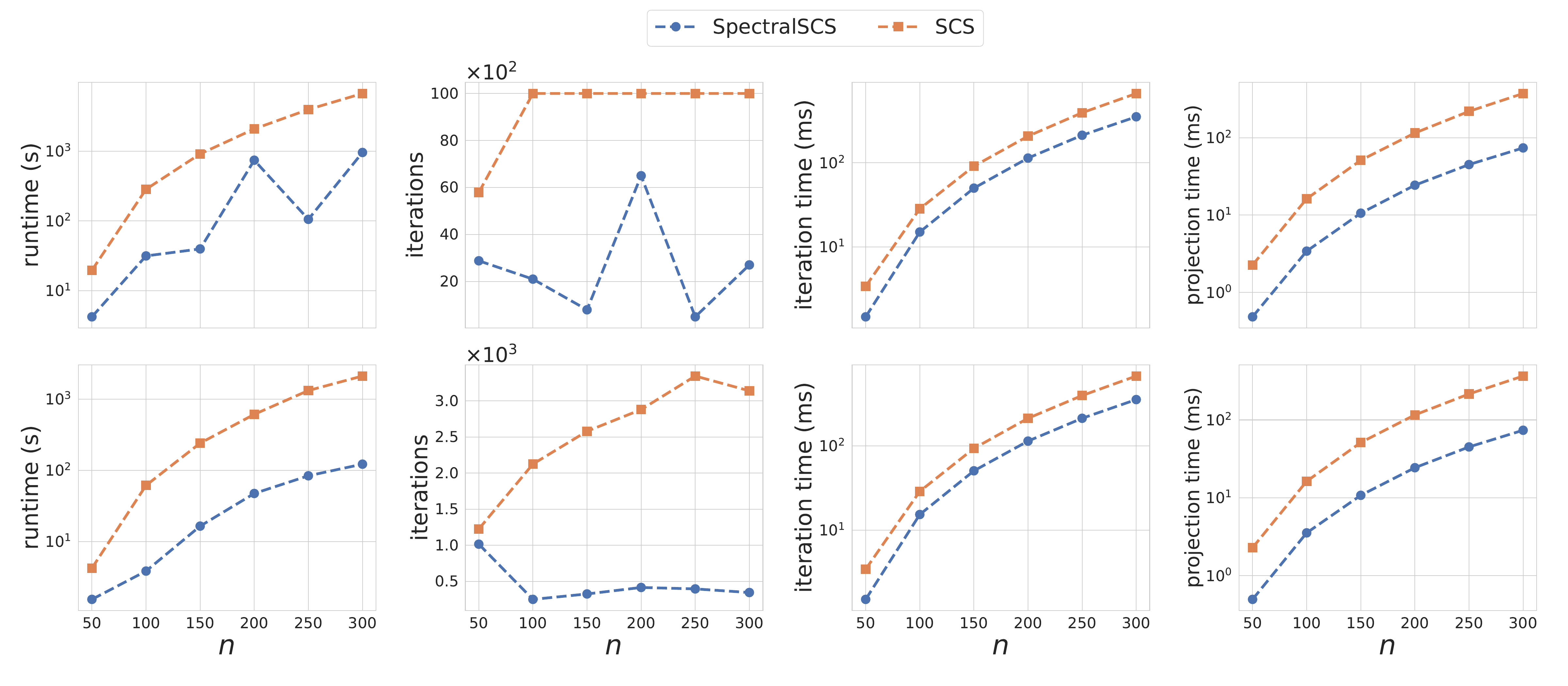}
\caption{Results for experimental design for standard accuracy 
$\epsilon_{\text{abs}} = \epsilon_{\text{rel}} = 10^{-4}$ (top row) and
lower accuracy $\epsilon_{\text{abs}} = \epsilon_{\text{rel}} = 10^{-3}$ (bottom
row). The first column shows the total solve time, the second column shows the
total number of iterations, the third column shows the time per iteration, and
the fourth column shows the time per matrix cone projection.}
\label{fig:exp-design}
\end{figure}

The termination criterion for \texttt{SCS} and \texttt{SpectralSCS} depends on
two tolerances $\epsilon_{\text{abs}} > 0$ and $\epsilon_{\text{rel}} > 0$. The
first row in Figure~\ref{fig:exp-design} shows the result for $n \in \{50, 100,
\dots, 300\}$ for the default tolerance $\epsilon_{\text{abs}} =
\epsilon_{\text{rel}} = 10^{-4}$. We noticed that \texttt{SCS} struggled to
achieve the default accuracy within a pre-specified maximum of $10^4$
iterations. We therefore also show the result for lower accuracy
$\epsilon_{\text{abs}} = \epsilon_{\text{rel}} = 10^{-3}$ in the second row of
Figure~\ref{fig:exp-design}. Both for the default and the lower accuracy,
\texttt{SpectralSCS} converges with significantly fewer iterations and is often
an order of magnitude faster than \texttt{SCS}. On average, \texttt{SpectralSCS}
is 20.7 times faster than \texttt{SCS}. Interestingly, while
\texttt{SpectralSCS} exhibits much faster cone projections compared to
\texttt{SCS}, this improvement does not fully translate into a proportional
reduction in time per iteration. For example, for $n = 300$, the spectral
matrix cone projection of \texttt{SpectralSCS} is six times faster than the
positive semidefinite cone projection of \texttt{SCS}, but the overall iteration
speed of \texttt{SpectralSCS} is only roughly twice as fast. The reason for
this discrepancy is that for \texttt{SpectralSCS}, the linear system solve with
cached factorization begins to dominate. 

\subsubsection{Sparse inverse covariance selection}
We consider the task of estimating a covariance matrix $\Sigma \in \symm^n$,
under the prior assumption that $\Sigma^{-1}$ is sparse. Given a sample
covariance matrix $S \in \symm^n_+$, this \emph{covariance selection problem} can
be formulated as \cite{Friedman2007}
\BEQ \label{e:sparse-inv-prob}
\begin{array}{ll}
\mbox{minimize} & \Tr(SX) - \log \det X + \lambda \|X \|_1,
\end{array} 
\EEQ
where the variable $X \in \symm^n$ is the estimate of $\Sigma^{-1}$, $\lambda >
0$ is a regularization parameter, and $\| X \|_1 = \sum_{i = 1}^n \sum_{j=1}^n
|X_{ij}|$. To solve \eqref{e:sparse-inv-prob} with \texttt{SpectralSCS} we
formulate it as 
\[
\begin{array}{ll}
\mbox{minimize} & \Tr(SX) + t + \lambda \sum_{i=1}^n \sum_{j \geq i}^n z_{ij} \\
\mbox{subject to} &  -z_{ii} \leq X_{ii} \leq z_{ii}, \hspace{0.6cm} i = 1, \dots, n \\
& -z_{ij} \leq 2 X_{ij} \leq z_{ij}, \hspace{0.3cm} i = 1, \dots, n, \: j > i \\
& v  = 1 \\
& (t, v, X) \in K_{\text{logdet}},
\end{array}, 
\]
with variables $t \in \reals$, $v \in \reals$, $X \in \symm^n$, and $z \in
\reals^{n(n+1)/2}$. The main cost of the matrix cone projection for \texttt{SCS}
and \texttt{SpectralSCS} is to compute the eigenvalue decomposition of a
symmetric matrix of dimension $2n$ and $n$, respectively.

To generate the problem data $S$ we use the same procedure as in
\cite{wang2010}. The true inverse covariance matrices have a density of
non-zero entries around $5\%$, and the regularization parameter $\lambda$ was
chosen to achieve similar density for the estimates.

\begin{figure}[!htb]
\centering
\includegraphics[width=1\textwidth]{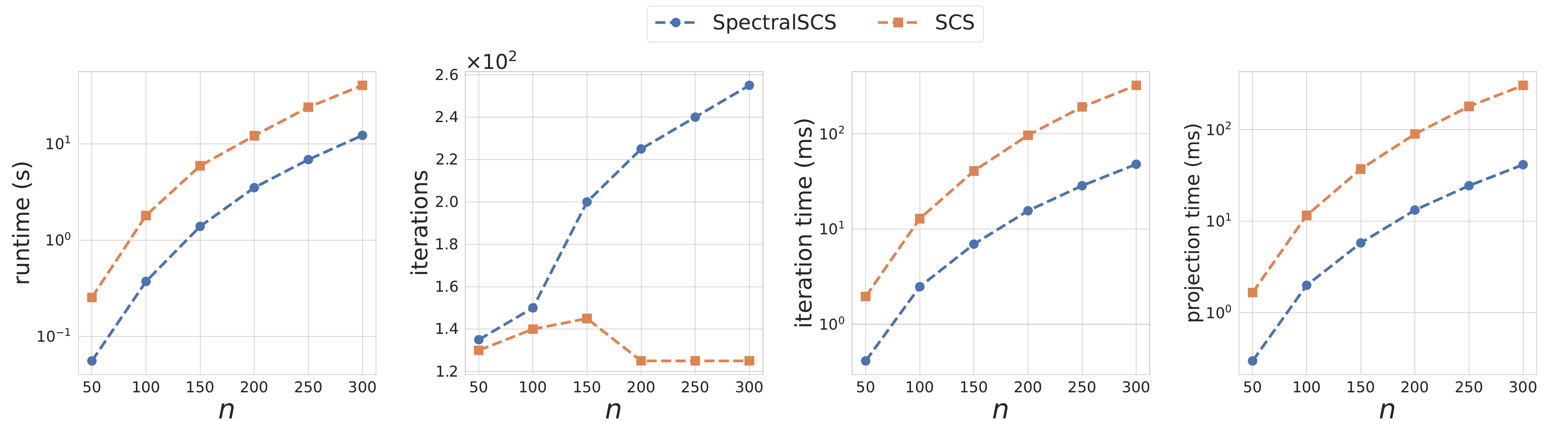}
\caption{Results for sparse inverse covariance selection.}
\label{fig:sparse_inv}
\end{figure}

Figure~\ref{fig:sparse_inv} shows the result for $n \in \{50, 100, \dots,
300\}$. The spectral matrix cone projection of \texttt{SpectralSCS} is
significantly faster than the positive semidefinite cone projection of
\texttt{SCS}, which translates to an equal reduction in time per iteration. For
example, for $n = 300$, the spectral matrix cone projection is six times faster,
and \texttt{SpectralSCS} has six times faster iterations than \texttt{SCS}.
However, \texttt{SpectralSCS} requires more iterations to converge, but on
average it is still 4.0 times faster than \texttt{SCS}.

\subsection{Nuclear norm cone}
We consider the problem of robust principal component analysis, as introduced in
\S \ref{sec:potential advantage}. The main cost of \texttt{SCS} in each iteration is to
compute the eigenvalue decomposition of a matrix of size $(m + n) \times (m +
n)$. For \texttt{SpectralSCS} the main cost is to compute the reduced SVD 
of a matrix of size $m \times n$. We show results for $m = n, \: m
= 2n$ and $m = 5n$, where we vary $n$. The problem instances are generated as in
\cite{ODonoghue2016}. Specifically, we set $M = \hat{L} + \hat{S}$ where $\hat{L}$ was a
randomly generated rank-$r$ matrix and $\hat{S}$ was a sparse matrix with
approximately 10\% nonzero entries. For all instances, we set $\mu$ to be equal
to the sum of absolute values of the entries of $\hat{S}$ and generated the data
with $r = 10$.

\begin{figure}[!htb]
\centering
\includegraphics[width=0.95\textwidth]{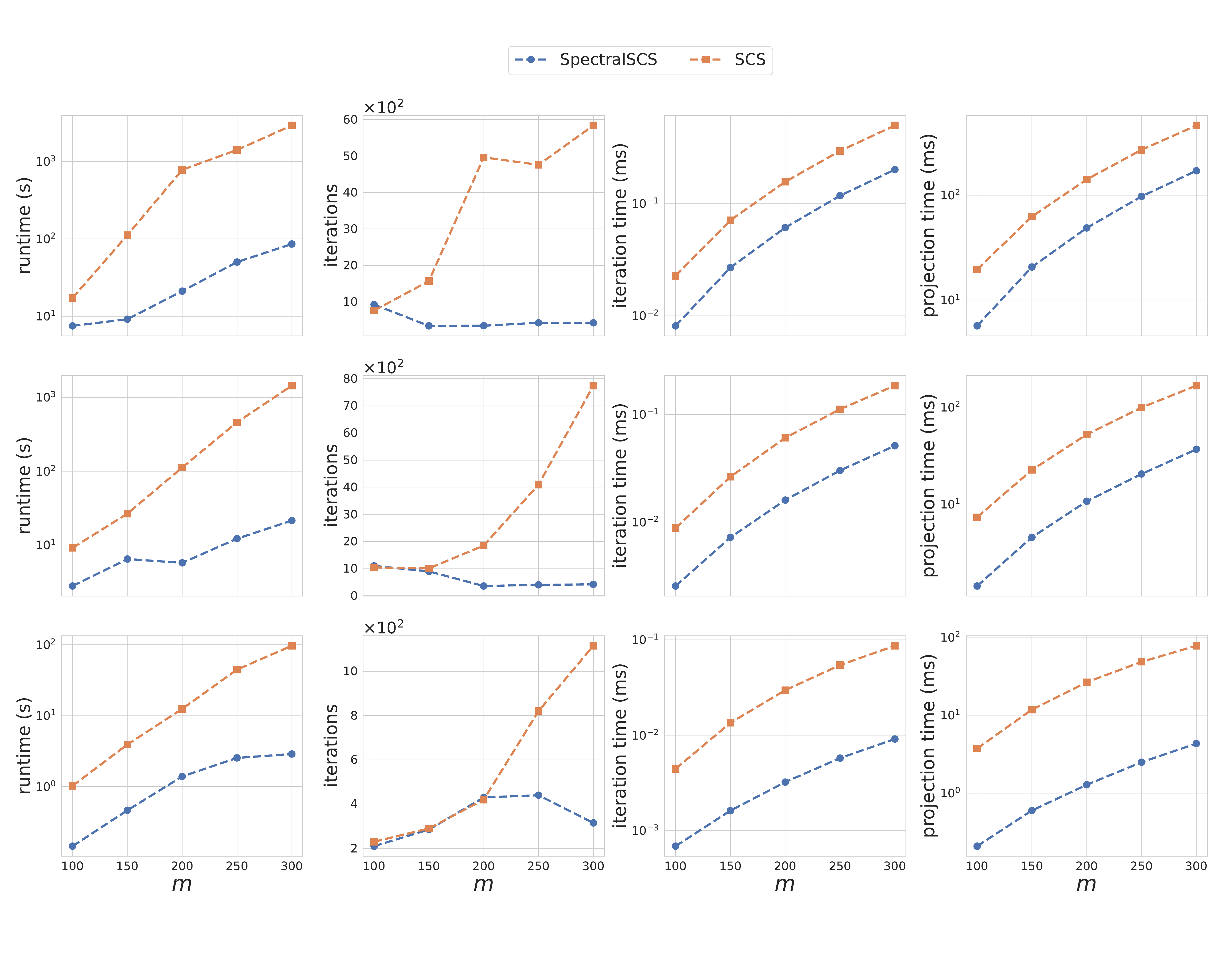} 
\caption{Results for robust PCA for $m = n$ (first row),
$m = 2n$ (second row), and $m = 5n$ (third row).
\label{fig:nuc_norm_main}}
\end{figure}

Figure~\ref{fig:nuc_norm_main} shows the results. We see that
\texttt{SpectralSCS} converges with significantly fewer iterations and is often
an order of magnitude faster than \texttt{SCS}. On average, \texttt{SpectralSCS}
is 22.3 times faster than \texttt{SCS}. For $m = 300$ and $m = n, \: m = 2n$,
and  $m = 5n$, the iterations of \texttt{SpectralSCS} are 2.6, 3.6 and 8.6 times
faster than the iterations of \texttt{SCS}, respectively.

\subsection{Sum-of-largest-eigenvalues cone}
We consider the following variant of the so-called \emph{graph partitioning
problem}. The problem is to divide $n$ nodes of a given graph into $k$ disjoint
subsets of equal size such that the total number of edges connecting different 
subsets is minimized. This problem is NP-hard, but Donath and Hoffman \cite{Donath1972}
have derived a procedure for obtaining a suboptimal partitioning using 
the eigenvectors of a matrix $\diag(x^\star) - L$, where $L$ is the Laplacian
of the graph and $x^\star$ is the solution of 
\BEQ \label{e:graph-p-problem}
\begin{array}{ll}
  \mbox{minimize} & \sum_{i=1}^k \boldsymbol{\lambda}_i(\diag(x) - L) \\
  \mbox{subject to} &  \ones^T x = 0,
\end{array}
\EEQ
with variable $x \in \reals^n$. To solve \eqref{e:graph-p-problem}
with \texttt{SpectralSCS} we formulate it as 
\[
\begin{array}{ll}
  \mbox{minimize} & t\\
  \mbox{subject to} &  \ones^T x = 0 \\
  & (t, \diag(x) - L) \in K_{\text{mSum}},
\end{array} 
\]
with variables $t \in \reals$ and $x \in \reals^n$. The main cost of the matrix
cone projection for \texttt{SCS} is to compute \emph{two} eigenvalue
decompositions of symmetric matrices of dimension $n$, while the main cost for
\texttt{SpectralSCS} is to compute \emph{one} eigenvalue decomposition of a
symmetric matrix of dimension $n$.

To generate the problem data $L$ we construct graphs with $n$ nodes using the Erdos-Renyi
random graph model with edge probability $p = 0.01$ \cite{erdos1959}.  

\begin{figure}[!htb]
\centering
\includegraphics[width=0.95\textwidth]{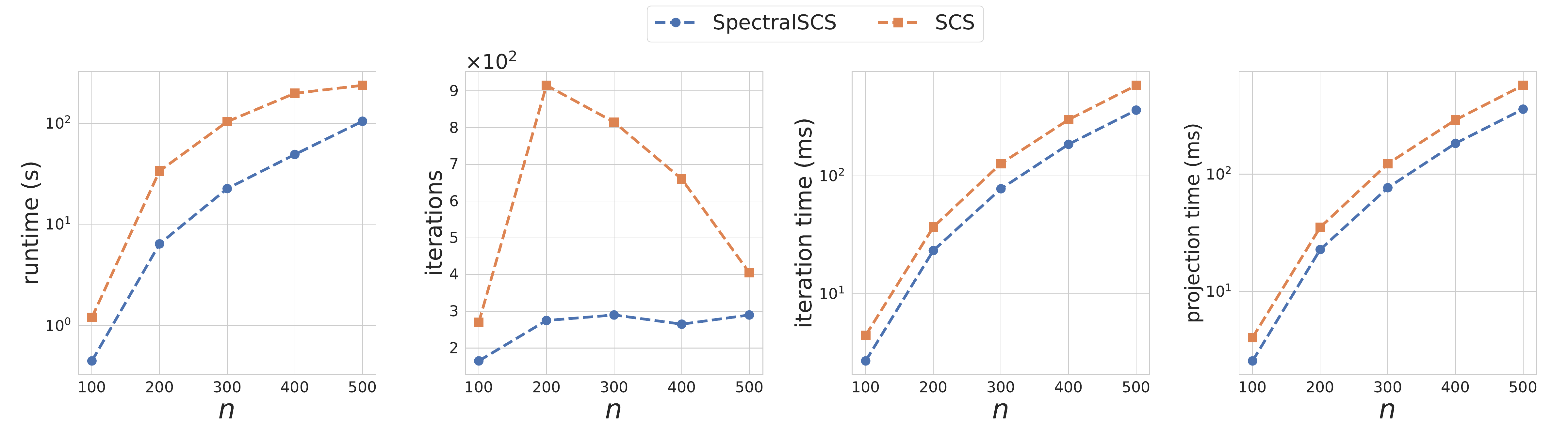}
\caption{Results for graph partitioning.}
\label{fig:graph_part}
\end{figure}

Figure~\ref{fig:graph_part} shows the result for $n \in \{100, 200, \dots,
500\}$ and $k = 10$. We see that \texttt{SpectralSCS} converges within fewer
iterations and has faster iterations, resulting in a total speedup compared to 
\texttt{SCS}. On average, \texttt{SpectralSCS} is 3.9 times faster than \texttt{SCS}.

\subsection{Overhead of spectral vector cone projection}
\label{sec:ablation}
\begin{figure}[!htb]
  \centering
  \subfloat[]{\includegraphics[width=0.31\textwidth]{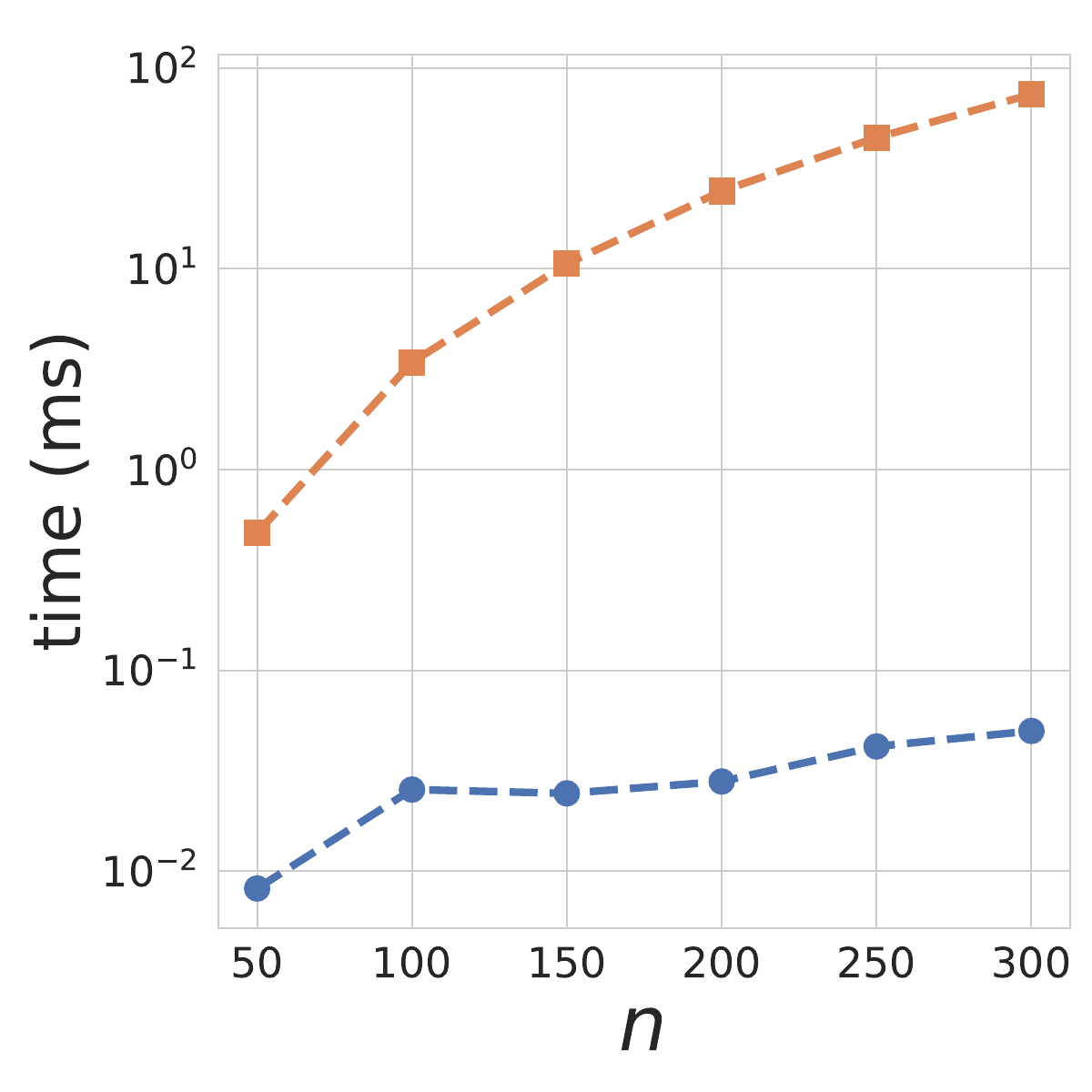}}
  \subfloat[]{\includegraphics[width=0.31\textwidth]{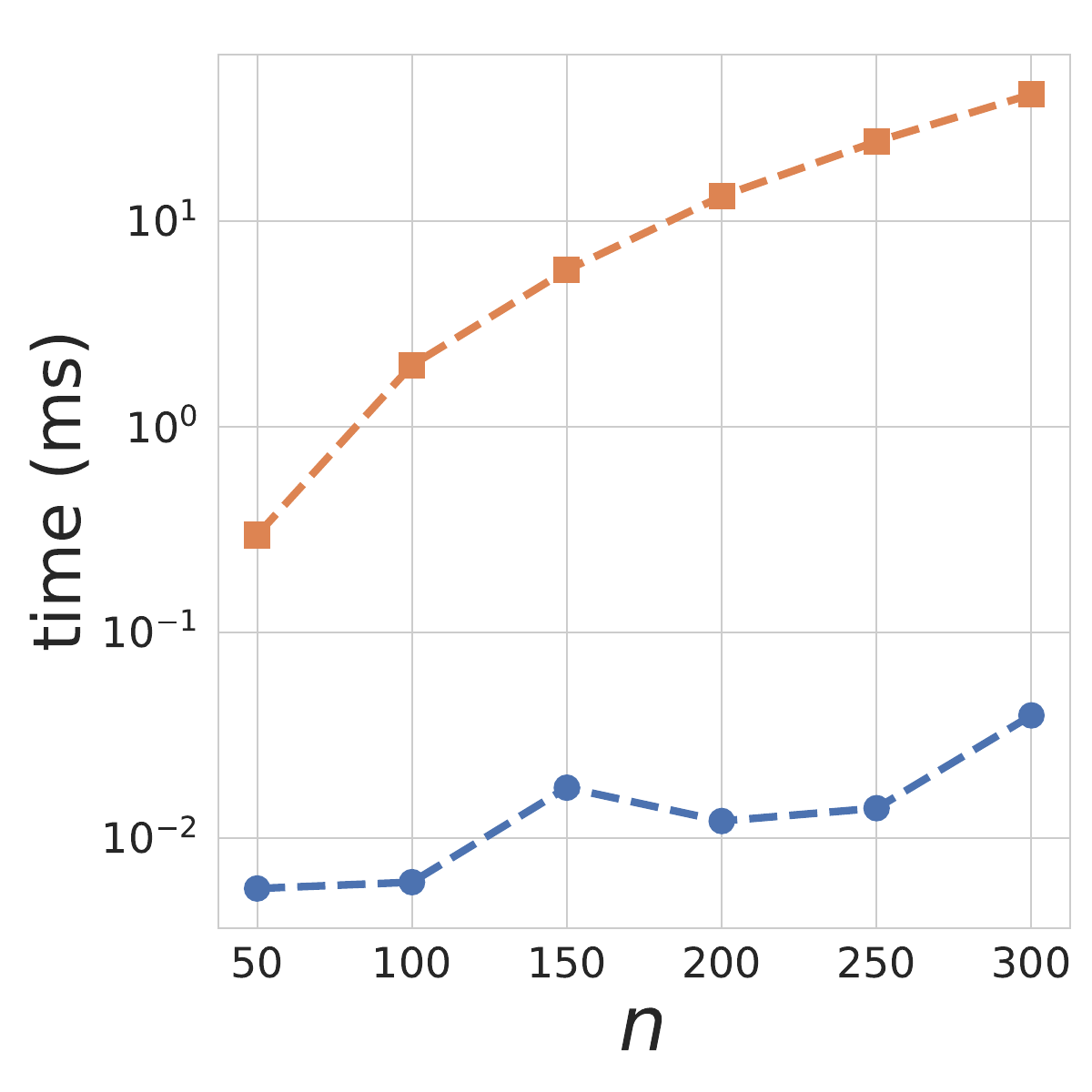}}
  \subfloat[]{\includegraphics[width=0.31\textwidth]{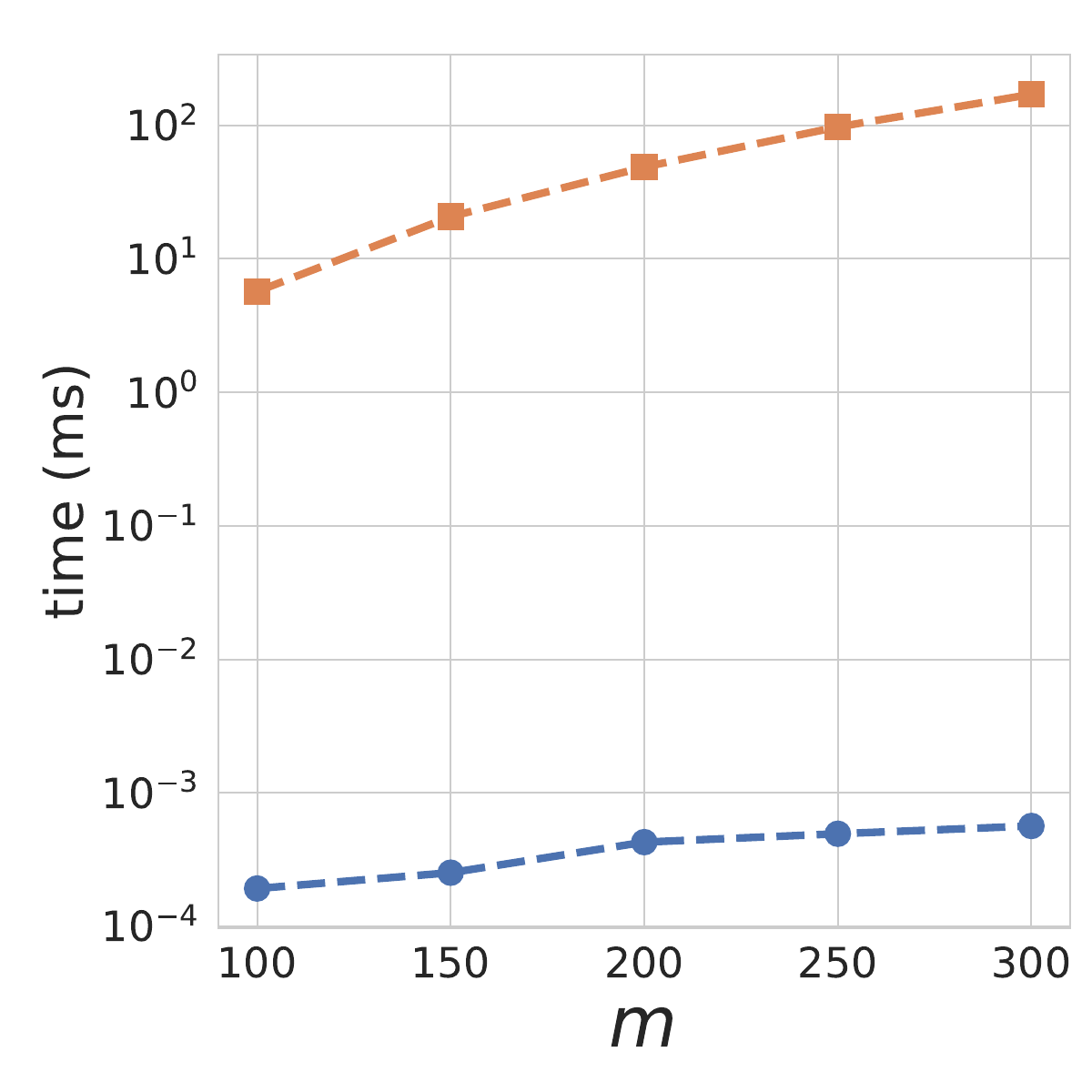}}
  \\
  \subfloat[]{\includegraphics[width=0.31\textwidth]{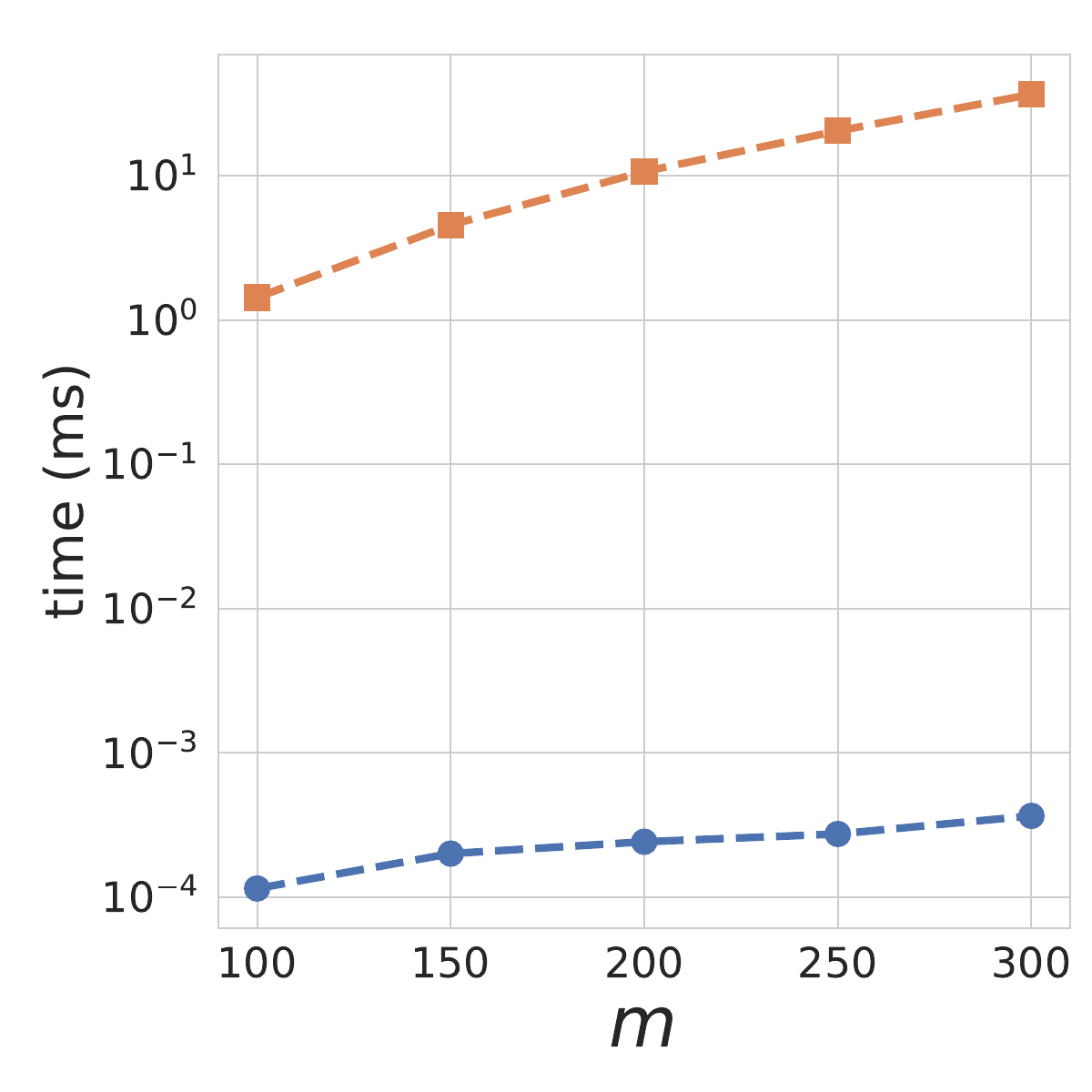}}
  \subfloat[]{\includegraphics[width=0.31\textwidth]{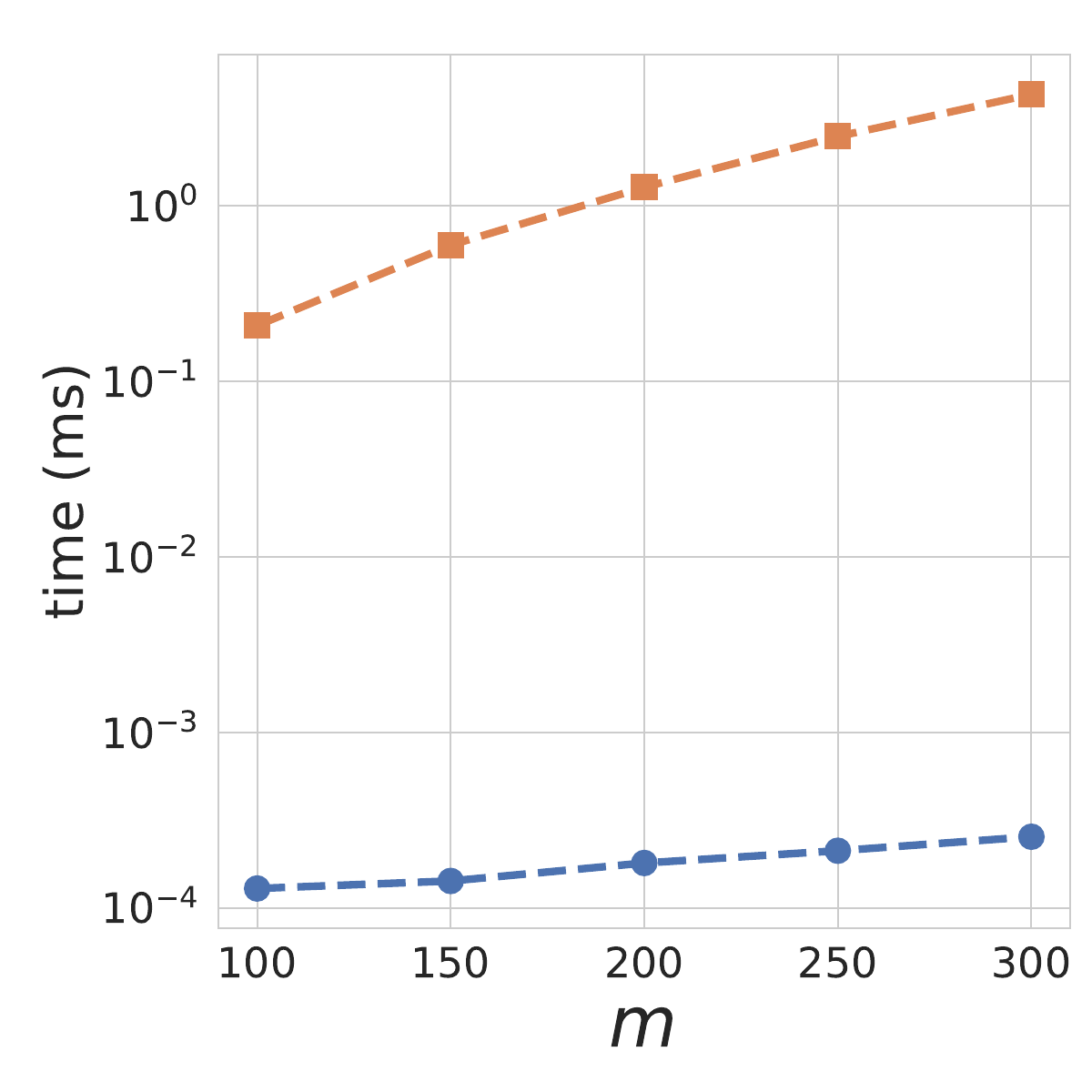}}
  \subfloat[]{\includegraphics[width=0.31\textwidth]{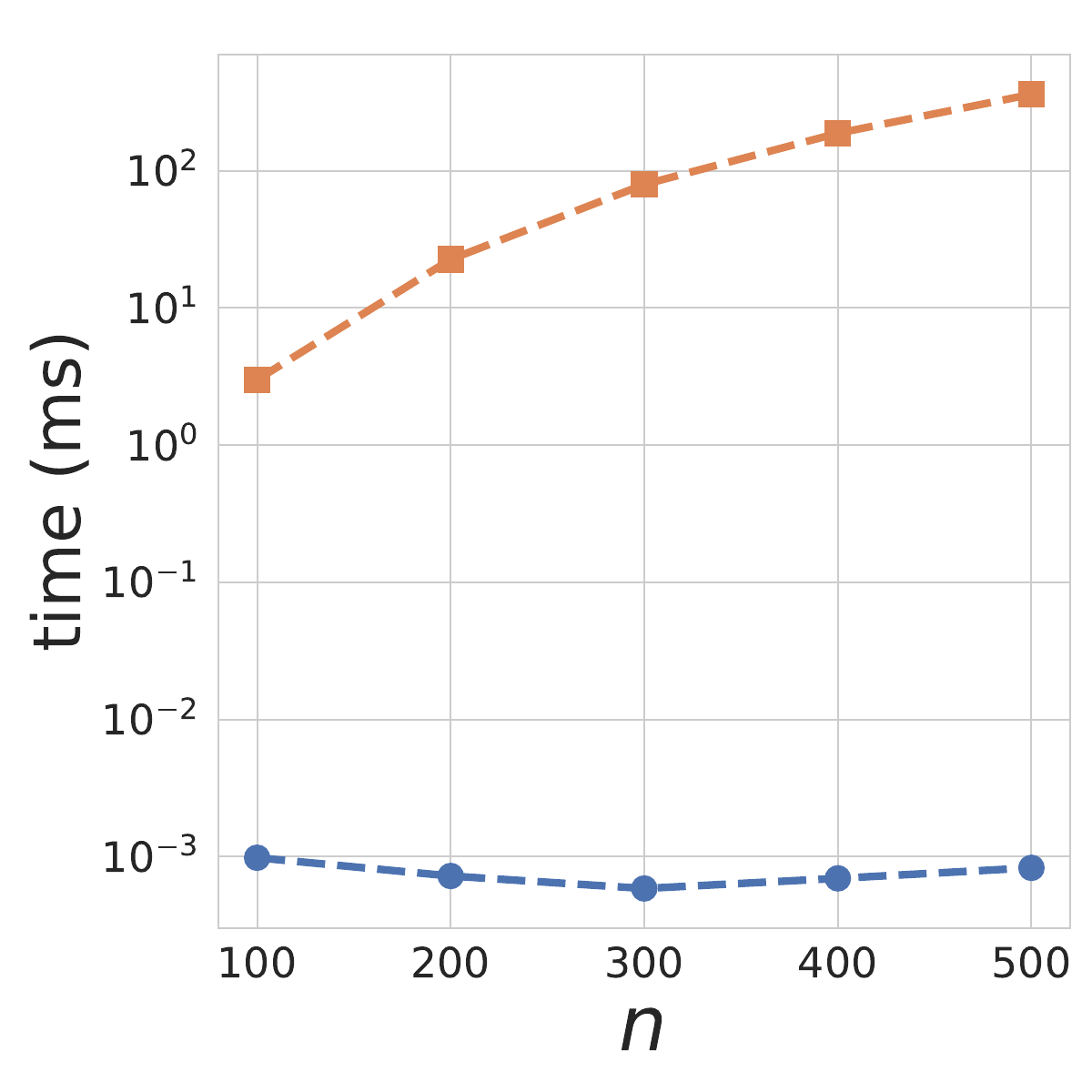}}
  \caption{The average time to compute the eigenvalue or singular value decomposition
  (orange line) and the average time to project onto the spectral vector cone (blue line)
  for \textbf{(a)} experimental design, \textbf{(b)} sparse inverse covariance
  selection, 
  \textbf{(c)} robust
  PCA with $m = n$, \textbf{(d)} robust PCA with $m = 2n$, \textbf{(e)} robust PCA with $m = 5n$,
  and \textbf{(f)} graph partitioning.}
  \label{fig:ablation}
  \end{figure}

For spectral matrix cones to be an attractive alternative to the standard
canonicalization based on the positive semidefinite cone, it is crucial that the
spectral \emph{vector} cone projection does not introduce a significant
overhead. In Figure~\ref{fig:ablation} we compare the average time to compute
the eigenvalue or singular value decomposition required for the spectral
\emph{matrix} cone projection with the average time for projecting the
eigenvalues or singular values onto the corresponding spectral \emph{vector}
cone. We see that the spectral vector cone projection is often at least two
orders of magnitude faster, making it insignificant in comparison.

\newpage
\section{Conclusions and extensions}
We considered the class of spectral matrix cones, and showed that projecting a
matrix onto such a cone can be done by projecting its eigenvalues or singular
values onto an associated spectral vector cone. We showed that these spectral
vector cone projections can be implemented efficiently, rendering this
projection step negligble when compared to the time required for the eigenvalue
or singular value decomposition. By integrating the new projection routines into
the first-order conic solver SCS, we observed a significant improvement in
solver performance.

\paragraph{Extensions.}
Spectral matrix cones can be used to canonicalize many interesting matrix
functions, such as the nuclear norm or the log-determinant function. Another new
matrix cone that could potentially improve first-order conic solvers like SCS is
the $\ell_1$-matrix cone $K = \{ (t, X) \in \reals \times \symm^n \: | \: \|X
\|_1 \leq t\}$. This is not a spectral matrix cone, but projecting onto it can
be done efficiently using the algorithm described in \S \ref{sec:ad-hoc
projection} for projecting onto $K_{\ell_1}$. Similarily, the projection onto the 
matrix box cone 
$K = \{(t, X) \in \reals \times \symm^n \: | \: t \ell I \preceq X \preceq t u I\}$
(which together with the constraint $t = 1$ can be used to enforce bounds on the 
eigenvalues of a matrix) can be computed efficiently based on ideas similar
to those in this paper.

Another matrix cone that has the potential to significantly improve SCS for
solving certain semidefinite programs is the cone of positive semidefinite
completable matrices with a given chordal sparsity pattern
\cite{Vandenberghe2015}. It is relatively expensive to compute the
\emph{Euclidean} projection onto this cone \cite{Sun2015}, but a
\emph{generalized projection} with respect to a carefully chosen Bregman
divergence can be computed at a low cost corresponding to a few sparse Cholesky
factorizations \cite{Jiang2022}.
\label{sec:conclusions}

Finally, an interesting direction for future work is extending the theory of 
differentiating the solution map of cone programs \cite{Agrawal2019} by
exploring whether and how the differential of spectral matrix cone projectors 
can be computed.

\newpage
\appendix 

\section{Explicit cone expressions}
\label{sec:appendix-topology}
In Appendix \ref{sec:appendix-topology} we explicitly characterize a few 
spectral vector cones and get rid of the closure 
in their definition. 

\subsection{Logarithmic cone}
\label{sec:appendix-topology-log} 
We prove that 
\[
\begin{split}
K_{\text{log}} & = \cl \{(t, v, x) \in \reals \times \reals_{++}
\times \reals_{++}^n \: | \: -\sum_{i=1}^n v \log (x_i/v) \leq t \} =
K_{\text{log}}^0 \cup (\reals_+ \times \{0\} \times \reals_+^n)
\end{split}  
\] 
where 
\[
K_{\text{log}}^0 = \{(t, v, x) \in \reals \times \reals_{++}
\times \reals_{++}^n \: | \: - \sum_{i=1}^n v \log (x_i/v) \leq t \}.
\]
Let $\overline{K}_{\text{log}} =K_{\text{log}}^0 \cup (\reals_+ \times \{0\}
\times \reals_+^n).$ 
We first prove that
$\overline{K}_{\text{log}} \subseteq K_{\text{log}}$. Consider a point 
$(\bar{t}, \bar{v}, \bar{x}) \in \reals_+ \times \{0\} \times \reals_+^n$.  
(If $(\bar{t}, \bar{v},
\bar{x}) \in K^0_{\text{log}}$ it clearly holds that $(\bar{t}, \bar{v},
\bar{x}) \in K_{\text{log}}$.) For $k \geq 1$, define the sequence $(t^{(k)}, v^{(k)},
x^{(k)})$
with 
\[
t^{(k)} = \bar{t}, \qquad v^{(k)} = 1/k, \qquad x^{(k)} =  \bar{x} + 1/k.    
\] Then $(t^{(k)}, v^{(k)}, x^{(k)}) \to (\bar{t}, \bar{v}, \bar{x})$ and 
$(t^{(k)}, v^{(k)}, x^{(k)}) \in K^0_{\text{log}}$ for all $k \geq 1$ since 
\[
-\sum_{i=1}^n v^{(k)} \log(x_i^{(k)}/v^{(k)}) = - \sum_{i=1}^n v^{(k)} 
\log \bigg(\frac{\bar{x}_i + 1/k}{1/k} \bigg)
\leq - \sum_{i=1}^n v^{(k)} \log\bigg( \frac{1/k}{1/k} \bigg) = 0 
\leq \bar{t} = t_k.    
\] 
Hence, $(\bar{t}, \bar{v}, \bar{x})$ is an accumulation point of a sequence 
in $K^0_{\text{log}}$, implying that $(\bar{t}, \bar{v}, \bar{x}) \in K_{\text{log}}.$

Conversely, let us show that $K_{\text{log}} \subseteq
\overline{K}_{\text{log}}$. Consider a point 
$(\bar{t}, \bar{v}, \bar{x}) \in K_{\text{log}}$ and assume that 
$(\bar{t}, \bar{v}, \bar{x}) \not \in \reals_+ \times \{0\} \times \reals_+^n$.
(If $(\bar{t}, \bar{v}, \bar{x}) \in \reals_+ \times \{0\} \times \reals_+^n$
we are done.)
We shall prove that $(\bar{t}, \bar{v}, \bar{x}) \in K^0_{\text{log}}$. 

If $(\bar{t}, \bar{v}, \bar{x}) \in K_{\text{log}}$ it must hold that $\bar{v} \geq
0$ and $\bar{x} \in \reals^n_+$. Hence, the only possible situation in which 
$(\bar{t}, \bar{v}, \bar{x}) \in K_{\text{log}}$ but
$(\bar{t}, \bar{v}, \bar{x}) \not \in \reals_+ \times \{0\} \times \reals_+^n$
can happen is if either $v > 0$ or $t < 0$.

Since $(\bar{t}, \bar{v}, \bar{x}) \in K_{\text{log}}$ there
exists a sequence $(t^{(k)}, v^{(k)}, x^{(k)}) \subset \reals \times \reals_{++}
\times \reals_{++}^n$ such that $(t^{(k)}, v^{(k)}, x^{(k)}) \to (\bar{t},
\bar{v}, \bar{x})$ and 
\[
- \sum_{i=1}^n v^{(k)} \log(x_i^{(k)}/v^{(k)}) \leq t^{(k)}.    
\]  
Assume $v > 0$. We must show that $x \in \reals^n_{++}$. Assume $x_j = 0$ for
some $j \in \{1, \dots, n\}$. Then the left side of the inequality converges to $\infty$
as $k \to \infty$ and the right side converges to the finite value $\bar{t}$.
Hence, $x \in \reals^{n}_{++}$ must hold, implying that $(\bar{t}, \bar{v}, \bar{x}) \in
K^0_{\text{log}}$. Now assume $\bar{t} < 0$. We want to show that $\bar{v} > 0$
and $\bar{x} \in \reals^{n}_{++}$. There are three cases:
\BIT 
\item If $\bar{x}_j = 0$ for some $j \in \{1, \dots, n\}$ and $\bar{v} > 0$, then the left side of the inequality converges to
$\infty$ and the right side converges to $\bar{t} < 0$.
\item If $\bar{x} \in \reals^n_{++}$ and $\bar{v} = 0$, then the left side of the inequality
converges to 0
and the right side converges to $\bar{t} < 0$.
\item If $\bar{x}_j = 0$ for some $j \in \{1, \dots, n\}$ and $\bar{v} = 0$, then the left side of the inequality converges to
something nonnegative. 
\EIT  In all cases we get a contradiction of the inequality, implying that $\bar{v} >
0$ and $\bar{x} \in \reals^n_{++}$ so $(\bar{t}, \bar{v}, \bar{x}) \in K^0_{\text{log}}$. 

\subsection{Inverse cone}
\label{sec:appendix-topology-inverse}
We prove that 
\[
\begin{split}
K_{\text{inv}} & = \cl \{(t, v, x) \in \reals \times \reals_{++}
\times \reals^n_{++} \: | \: v^2 \sum_{i=1}^n 1/x_i \leq t \} = 
K_{\text{inv}}^0 \cup (\reals_+ \times \{ 0 \} \times \reals_{+}^n) 
\end{split}  
\] 
where 
\[
K^0_{\text{inv}} = \{(t, v, x) \in \reals \times \reals_{++}
\times \reals^n_{++} \: | \: v^2 \sum_{i=1}^n 1/x_i \leq t \}.   
\]
Let 
$\overline{K}_{\text{inv}} =
K_{\text{inv}}^0 \cup (\reals_+ \times \{ 0 \} \times \reals_{+}^n) .$  
We first prove that
$\overline{K}_{\text{inv}} \subseteq K_{\text{inv}}$. Consider a point 
$(\bar{t}, \bar{v}, \bar{x}) \in \reals_+ \times \{0 \} \times \reals^n_+.$
(If $(\bar{t}, \bar{v}, \bar{x}) \in K^0_{\text{inv}}$ it clearly holds that
$(\bar{t}, \bar{v}, \bar{x}) \in K_{\text{inv}}$.)
For $k \geq 1$, define the sequence $(t^{(k)}, v^{(k)}, x^{(k)})$ with 
\[
t^{(k)} = \bar{t} + 1/k, \qquad v^{(k)} = 1/k, \qquad x^{(k)} = \bar{x} + (n/k)\ones.    
\] 
Then $(t^{(k)}, v^{(k)}, x^{(k)}) \to (\bar{t}, \bar{v}, \bar{x})$ and 
$(t^{(k)}, v^{(k)}, x^{(k)}) \in K^0_{\text{inv}}$ for $k \geq 1$ since 
\[
(v^{(k)})^2 \sum_{i=1}^n \frac{1}{x^{(k)}_i} \leq \frac{1/k^2}{1/k} = \frac{1}{k} \leq t^{(k)}.
\] 
Hence, $(\bar{t}, \bar{v}, \bar{x})$ is an accumulation point of a sequence in $K^0_{\text{inv}}$,
implying that $(\bar{t}, \bar{v}, \bar{x}) \in K_{\text{inv}}.$

Conversely, let us show that $K_{\text{inv}} \subseteq
\overline{K}_{\text{inv}}$. Consider a point 
$(\bar{t}, \bar{v}, \bar{x}) \in K_{\text{inv}}$ and assume that 
$(\bar{t}, \bar{v}, \bar{x}) \not \in \reals_+ \times \{0\} \times \reals^n_+$. 
(If $(\bar{t}, \bar{v}, \bar{x}) \in \reals_+ \times \{0\} \times \reals^n_+$
we are done.)
We shall prove that $(\bar{t}, \bar{v}, \bar{x}) \in K^0_{\text{inv}}$. 
If $(\bar{t}, \bar{v}, \bar{x}) \in K_{\text{inv}}$ it must hold that $\bar{t} \geq 0, \bar{v} \geq 0$ and 
$\bar{x} \in \reals^n_+$. Hence, the only possible situation in which 
$(\bar{t}, \bar{v}, \bar{x}) \in K_{\text{inv}}$ but $(\bar{t}, \bar{v},
\bar{x}) \not \in \reals_+ \times\{0\} \times \reals^n_+$ can happen is if 
$\bar{t} \geq 0, \bar{v} > 0$ and $\bar{x} \in \reals^n_+$. There are two cases:
if $\bar{x} \in \reals^{n}_{++}$, or if $\bar{x}_j = 0$ for some $j \in \{1, \dots,
n\}$. If $\bar{x} \in
\reals^{n}_{++}$ it means that $(\bar{t}, \bar{v}, \bar{x}) \in
K_{\text{inv}}^0$ and we are done. Assume $\bar{x}_j = 0$ for some $j \in \{1,
\dots, n\}$. Since $(\bar{t}, \bar{v}, \bar{x}) \in K_{\text{inv}}$ there exists 
a sequence $(t^{(k)}, v^{(k)}, x^{(k)}) \subset K_{\text{inv}}^0$ such that 
$(t^{(k)}, v^{(k)}, x^{(k)}) \to (\bar{t}, \bar{v}, \bar{x})$. In particular,
this sequence satisfies 
\[
(v^{(k)})^2 \sum_{i=1}^n \frac{1}{x^{(k)}_i} \leq t^{(k)} 
\iff (v^{(k)})^2 \bigg(1 + x_j^{(k)} \sum_{i\neq j} \frac{1}{x^{(k)}_i} \bigg) \leq t^{(k)} x^{(k)}_j.
\]
If $k \to \infty$, the left side of the inequality converges to a strictly positive value,
whereas the right side of the inequality converges to zero, thus resulting in a
contradiction
of the assumption that $\bar{x}_j = 0$ for some $j \in \{1, \dots, n\}$. 

\subsection{Entropy cone}
\label{sec:appendix-topology-entropy}
We prove that
\[
\begin{split}
K_{\text{vEnt}} & = \cl \{(t, v, x) \in \reals \times \reals_{++}
\times \reals_{+}^n \: | \: \sum_{i=1}^n x_i \log(x_i/v_i) \leq t \} 
 = K_{\text{vEnt}}^0 \cup (\reals_+ \times \reals_+ \times \{0 \}^n)
\end{split}  
\] 
where 
\[
K_{\text{vEnt}}^0 = \{(t, v, x) \in \reals \times \reals_{++}
\times \reals_{+}^n \: | \: \sum_{i=1}^n x_i \log(x_i/v_i) \leq t \}.
\]
Let $\overline{K}_{\text{vEnt}} =  K_{\text{vEnt}}^0 \cup (\reals_+ \times
\reals_+ \times \{0 \}^n).$  We first prove that
$\overline{K}_{\text{vEnv}} \subseteq K_{\text{vEnt}}$. Consider a point 
$(\bar{t}, \bar{v}, \bar{x}) \in \reals_+ \times \reals_+ \times \{0\}^n$. (If $(\bar{t}, \bar{v},
\bar{x}) \in K^0_{\text{vEnt}}$ it clearly holds that $(\bar{t}, \bar{v},
\bar{x}) \in K_{\text{vEnt}}$.) For $k \geq 1$, define the sequence $(t^{(k)}, v^{(k)},
x^{(k)})$
with 
\[
t^{(k)} = \bar{t}, \qquad v^{(k)} = \bar{v} + 1/k, \qquad x^{(k)} = 0.    
\] Then $(t^{(k)}, v^{(k)}, x^{(k)}) \to (\bar{t}, \bar{v}, \bar{x})$ and 
$(t^{(k)}, v^{(k)}, x^{(k)}) \in K^0_{\text{vEnt}}$ for $k \geq 1$ since 
\[
\sum_{i=1}^n x^{(k)}_i \log(x^{(k)}_i/v^{(k)}) = 0 \leq  \bar{t} = t^{(k)}.  
\]
Hence, $(\bar{t}, \bar{v}, \bar{x})$ is an accumulation point of a sequence in $K^0_{\text{vEnt}}$,
implying that $(\bar{t}, \bar{v}, \bar{x}) \in K_{\text{vEnt}}.$

Conversely, let us show that $K_{\text{vEnt}} \subseteq
\overline{K}_{\text{vEnt}}$. Consider a point 
$(\bar{t}, \bar{v}, \bar{x}) \in K_{\text{vEnt}}$ and assume that 
$(\bar{t}, \bar{v}, \bar{x}) \not \in \reals_+ \times \reals_+ \times \{0\}^n$.
(If $(\bar{t}, \bar{v}, \bar{x}) \in \reals_+ \times \reals_+ \times \{0\}^n$
we are done.), 
We shall prove that $(\bar{t}, \bar{v}, \bar{x}) \in K^0_{\text{vEnt}}$. 
If $(\bar{t}, \bar{v}, \bar{x}) \in K_{\text{vEnt}}$ it must hold that $\bar{v} \geq
0$ and $\bar{x} \in \reals^n_+$. Hence, the only possible situation in which 
$(\bar{t}, \bar{v}, \bar{x}) \in K_{\text{vEnt}}$ but
$(\bar{t}, \bar{v}, \bar{x}) \not \in \reals_+ \times \reals_+ \times \{0\}^n$
can happen is if either $\bar{x}_j > 0$ for some $j \in \{1, \dots, n\}$ or if $\bar{t} < 0$.

Since $(\bar{t}, \bar{v}, \bar{x}) \in K_{\text{vEnt}}$
there exists a sequence $(t^{(k)}, v^{(k)}, x^{(k)}) \subset K^0_{\text{vEnt}}$ such
that $(t^{(k)}, v^{(k)}, x^{(k)}) \to (\bar{t}, \bar{v}, \bar{x})$. In particular, this sequence satisfies 
\[
\sum_{i=1}^n x_i^{(k)} \log (x_i^{(k)}/v^{(k)}) \leq t^{(k)}.    
\] 
First assume $\bar{x}_j > 0$ for some $j \in \{1, \dots, n\}$. If $\bar{v} = 0$
the left side of the inequality converges to $\infty$ as $k \to \infty$, whereas
the right side converges to a finite value. Hence, it must hold that $\bar{v} >
0$, implying that $(\bar{t}, \bar{v}, \bar{x}) \in K^0_{\text{vEnt}}$. Now
assume $\bar{t} < 0$. We want to show that $\bar{v} > 0$. Assume $\bar{v} = 0$. 
If $\bar{x}_j > 0$ for some $j \in \{1, \dots, n\}$
and we let $k \to \infty$, the left side converges to $\infty$ and the right
side to $\bar{t} <0$. If $\bar{x} = 0$ the left side is 0 and the right side 
converges to $\bar{t} <0$. In both cases we get a
contradiction of the inequality, implying that $\bar{v} > 0$ so $(\bar{t}, \bar{v},
\bar{x}) \in K^0_{\text{vEnt}}$. 

\section{Dual cones}
\label{sec:appendix dual cones}
In Appendix \ref{sec:appendix dual cones} we derive dual cones of the spectral 
vector and spectral matrix cones presented in the main text. The derivations are 
based on the following identities from \S\ref{sec:spectral matrix cone
projections}:
\[
\begin{split}
K^*_f & =  \cl \{(t, v, x) \in \reals_{++} \times
\reals \times \reals^n \: | \: v \geq t f^*(-x/t)\} \\
K^*_F & =  \cl \{(t, v, X) \in \reals_{++} \times
\reals \times \symm^n \: | \: v \geq t F^*(-X/t)\} \\
F^*(X) & = f^*(\boldsymbol{\lambda}(X)).
\end{split}
\]

\paragraph{Log-determinant cone.} The conjugate of $f(x) = - \sum_{i=1}^n \log
x_i$, $\dom f = \reals^{n}_{++}$ is $f^*(y) = - n - \sum_{i=1}^n \log(-y_i)$, 
$\dom f^* = - \reals^{n}_{++}$. The dual of the logarithmic cone is therefore 
\[
K^*_{\text{log}} = \cl \{(t, v, x) \in \reals_{++} \times \reals \times \reals^{n}_{++}
\: | \: v \geq t(-n - \sum_{i=1}^n \log(x_i/t))\}.
\]
It can be expressed as 
\[
K^*_{\text{log}} = \{(t, v, x) \in \reals_{++} \times \reals \times \reals^{n}_{++}
\: | \: v \geq t(-n - \sum_{i=1}^n \log(x_i/t))\} \cup (\{0\} \times \reals_+ \times \reals^n_+).
\]
The conjugate of $F(X) = - \log \det X$, $\dom F = \symm^n_{++}$ is $F^*(Y) = - n - \log \det(-Y)$, 
$\dom F^* = - \symm^n_{++}$. The dual of the log-determinant cone is
therefore 
\[
K^*_{\text{logdet}} =
\cl\{(t, v, X) \in \reals_{++} \times \reals \times \symm^n_{++} \: | \: v \geq t(-n - \log \det(X/t))\}.
\]

\paragraph{Nuclear norm cone.} The conjugate of $f(x) = \| x \|_1$, $\dom f =
\reals^n$ is 
\[
f^*(y) = 
\begin{cases}
0 & \text{ if } \| x \|_\infty \leq 1 \\
\infty & \text{ otherwise}.
\end{cases}
\]
The dual of the $\ell_1$-norm cone is therefore 
\[
K^*_{\ell_1} = \{(t, x) \in \reals \times \reals^n \: | \: \| x \|_\infty \leq t \}.
\] 
(In general, the dual cone of a norm cone is the dual norm cone.)
The conjugate of $F(X) = \|X \|_*$, $\dom F = \reals^{m \times n}$ is 
\[
F^*(Y) = 
\begin{cases}
0 & \text{ if } \| Y \|_{2} \leq 1 \\
\infty & \text{ otherwise},
\end{cases}
\] 
where $\| Y \|_{2} = \sigma_1(Y)$ is the spectral norm.
The dual of the nuclear norm cone is therefore
\[
K_{\text{nuc}}^* = \{(t, X) \in \reals \times \reals^{m \times n} \: | \: \| X \|_2 \leq t\}.    
\]

\paragraph{Trace-inverse cone.} The conjugate of $f(x) = \sum_{i=1}^n 1/x_i$, 
$\dom f = \reals^{n}_{++}$ is $f^*(y) = -2\sum_{i=1}^n \sqrt{-y_i}$, $\dom f^* =
-\reals^n_+$.
The dual of the inverse cone is therefore 
\[
K_{\text{inv}}^* = \cl \{(t, v, x) \in \reals_{++} \times \reals \times  \reals^n_+
\: | \: v \geq - 2 t \sum_{i=1}^n \sqrt{x_i/t} \}.    
\]
The conjugate of $F(X) = \Tr(X^{-1})$, $\dom F = \symm^{n}_{++}$ is $F^*(Y) = -2 \Tr((-Y)^{1/2})$, 
$\dom F^* = - \symm^n_+.$ The dual of the trace-inverse cone is therefore 
\[
K^*_{\text{TrInv}} = \cl \{ (t, v, X) \in \reals_{++} \times \reals \times \symm^n_{+} \: | \:\
v \geq -2 \sqrt{t} \Tr(X^{1/2}) \}.
\]

\paragraph{Entropy cone.} The conjugate of $f(x) = \sum_{i=1}^n x_i \log x_i$,
$\dom f = \reals^n_+$ is $f^*(y) = \sum_{i=1}^n e^{y_i - 1}$, $\dom f^* =
\reals^n$. The dual of the vector entropy cone is therefore 
\[
K^*_{\text{vEnt}} = \cl \{ (t, v, x) \in \reals_{++} \times \reals \times \reals^n \: | \: 
v \geq t \sum_{i=1}^n e^{-x_i/t - 1}\}.    
\] 
The conjugate of $F(X) = \sum_{i=1}^n \lambda_i(X) \log \lambda_i(X), \: \dom F = \symm^n_+$
is $F^*(Y) = \Tr(\exp(Y - I))$, $\dom F^* = \symm^n$. The dual of the 
matrix entropy cone is therefore 
\[
K^*_{\text{mEnt}} = \cl \{ (t, v, X) \in \reals_{++} \times \reals \times \symm^n 
\: | \: v \geq t \Tr(\exp(-X/t - I))\}.
\] 
(Here $\exp(\cdot)$ is the matrix exponential.)

\paragraph{Root-determinant cone.} The conjugate of $f(x) = - \prod_{i=1}^n
x_i^{1/n}$, $\dom f = \reals^n_+$ is
\[
f^*(y) = 
\begin{cases}
0 & \text{ if } y \leq 0, \: \prod_{i=1}^n (-y_i)^{1/n} \geq 1/n \\
\infty & \text{ otherwise}.
\end{cases}
\]
The dual of the geometric mean cone is therefore 
\[
K^*_{\text{geomean}} = 
\cl\{ (t, x) \in \reals_{++} \times \reals^n_+ \: | \: (1/t) \prod_{i=1}^n x_i^{1/n} \geq 1/n\}.
\] 
The conjugate of $F(X) = - (\det(X))^{1/n}$, $\dom F = \symm^n_+$ is 
\[
F^*(Y) = 
\begin{cases}
0 & \text{ if } Y \preceq 0, \: (\det(-Y))^{1/n} \geq 1/n \\
\infty & \text{ otherwise}.
\end{cases}
\] 
The dual of the root-determinant cone is therefore 
\[
K_{\text{det}}^* = \cl \{(t, X) \in \reals_{++} \times \symm^n_+ \: | \: 
(1/t) (\det(X))^{1/n} \geq 1/n \}.
\]

\paragraph{Sum-of-largest-entries cone.} The conjugate of $f(x) = \sum_{i=1}^k
x_{[i]}$, $\dom f = \reals^n$ is 
\[
f^*(y) = 
\begin{cases}
0 & \text{ if } 0 \leq y \leq 1, \: \ones^T y = k \\
\infty & \text{ otherwise}.
\end{cases} 
\] 
The dual of the sum-of-largest-entries cone is therefore 
\[
K^*_{\text{vSum}} =  \{ (t, x) \in \reals \times \reals^n \: | \:
0 \geq x \geq -t, \: \ones^T x = -tk \}.
\] 
The conjugate of $F(X) = \sum_{i=1}^k \boldsymbol{\lambda}_i(X)$,
$\dom F = \symm^n$
is 
\[
F^*(Y) = 
\begin{cases}
0 & \text{ if } 0 \leq \boldsymbol{\lambda}(Y) \leq 1, \: 
\ones^T \boldsymbol{\lambda}(Y) = k \\
\infty & \text{ otherwise}.
\end{cases} 
\] 
The dual of the sum-of-largest-eigenvalues cone is therefore 
\[
K^*_{\text{mSum}} =  \{ (t, X) \in \reals \times \symm^n \: | \: 
0 \geq \boldsymbol{\lambda}(X) \geq -t, \: \ones^T \boldsymbol{\lambda}(X) = -tk \}.
\]

\section{Projecting onto the sum-of-largest-entries cone}
\label{sec:appendix sum-of-largest}
In Appendix \ref{sec:appendix sum-of-largest} we derive Algorithm 
\ref{alg:sum-largest} in \S \ref{sec:spectral vector cone
projections} for projecting onto the sum-of-largest-entries cone. The derivation
is inspired by \cite{Luxenberg25}, where the authors present an algorithm for
projecting onto the sublevel set of the sum-of-$k$-largest function.

The goal is to project $(\bar{t}, \bar{x})$ onto $K_{\text{vSum}}$, where
$\bar{x}$ is sorted with $\bar{x}_1 \geq \bar{x}_2 \geq \dots \geq \bar{x}_n$.
The idea behind the algorithm is to partition the entries of $\bar{x}$ into three
different blocks. The first block consists of \emph{untied entries}. Roughly
speaking, these are entries that when decreased by a sufficiently small amount
(while keeping all other entries fixed), reduce the sum of the $k$ largest
entries. The second block, which we call \emph{tied entries}, are the entries
that when decreased by a sufficiently small amount (while keeping all other
entries fixed), do not lead to a reduction in the sum of the $k$ largest entries
but instead reduce the $(k+1)$th largest unique entry of the vector. The last block
consists of the remaining entries.

We will maintain the number of entries in the
first two blocks, and the values of the tied entries and smallest untied entry,
using the following state variables:
\[
\begin{split}
n_u: \:  & \text{number of untied entries}, \hspace{1.7cm} n_t:
\text{number of tied entries} \\
a_u: \: & \text{value of smallest untied entry}, \qquad a_t: 
\text{value of tied entries}.
\end{split}
\]
We will also use a state variable $S$ to maintain the sum of the $k$ largest entries.

The projection algorithm is iterative. In most iterations, the algorithm decreases the tied
entries by a quantity $s > 0$ that depends on the current iteration. Simultaneously, it
decreases the untied entries by a quantity $sr$ and increases $t$ by $sr$, where
$r > 1$ is a parameter determined at each iteration (we will describe how to
select $r$ later).

The value of $s$ is chosen to ensure that one of two conditions is met: either
the cone constraint becomes satisfied, or a new tie occurs. The computation of
$s$ relies on the current state. The three scenarios below correspond to the
iteration (1) causing a tie between untied and tied entries, (2) causing a tie
between the tied entries and the largest entry of the last block, and (3)
making the cone constraint satisfied. Each scenario imposes a restriction on $s$, and we
then let $s = \min(s_1, s_2, s_3)$ where $s_1, \: s_2,$ and $s_3$ are the
restrictions from the three scenarios.
\BNUM 
\item If $n_u = k$, then the current iteration will decrease the untied entries
and leave the other entries unchanged until a tie occurs. This happens when
$a_u - s_1 = a_t$, so $s_1 = a_u - a_t$.  If there are no untied entries there is
no restriction from this scenario, so $s_1 = \infty$. If $1 \leq n_u \leq k-1$
we will decrease both the untied and tied entries, and a new tie between the
untied and tied elements occurs whenever
$a_t - s_1 = a_u - s_1 r$, so $s_1 = (a_u - a_t)/(r-1)$.
\item If there are no tied entries or if the last block is empty, then there is
no restriction from this scenario, so $s_2 = \infty$. Otherwise a tie between the
tied entries and the first element in the last block occurs when 
$a_t - s_2 = \bar{x}_{n_u + n_t + 1}$, so $s_2 = a_t - \bar{x}_{n_u + n_t + 1}$.
\item When the untied entries decrease by the quantity $s_3 r$, the sum of
the largest $k$ entries decreases by $n_u s_3 r$. When the tied entries decrease
by the quantity $s_3$, the sum of the largest $k$ entries decreases by $(k -
n_u)s_3$. Assuming that the vector remains sorted and no new ties occur, the
cone constraint is satisfied when 
\[
S - n_u s_3 r - (k - n_u) s_3 = t + s_3 r,
\] so $s_3 = (S-t)/(r(n_u + 1) + k - n_u)$.
\ENUM 
We should now decrease the untied entries by the quantity $sr$ where $s =
\min(s_1, s_2, s_3)$. To avoid vector subtraction in each iteration, we
introduce a new state variable,
\[
\eta: \text{the decrease to the untied entries}.
\]
In each iteration we update $\eta$ according to $\eta \gets \eta + s r$.
Decreasing the untied entries by $sr$ causes the sum of the $k$ largest entries
to decrease by $n_u s r$. Furthermore, when the tied entries decrease
by the quantity $s$, the sum of the largest $k$ entries decreases by $(k -
n_u)s$. We therefore update the state variable $S$ according to $S \gets S -
s(n_u r + k - n_u)$.

If there are tied entries (\ie~if $n_t > 0$), then we decrease them by $s$  so 
we let $a_t \gets a_t - s$. If a new tie between the untied and tied entries
occurred (\ie~if $s = s_1$), there is one less untied element so we let 
$n_u \gets n_u - 1$. If there is an untied element (\ie~if $n_u > 0$), we
conceptually decrease it by $\eta$ so we update $a_u \gets \bar{x}_{n_u} -
\eta$. Finally, we must update the number of tied entries. If there are no tied
entries (\ie~if $n_t = 0$) and there is a tie, there are two tied elements after
the tie so we let $n_t = 2$. If there is a tied element and a tie occurs, the
number of tied elements increases by one so we let $n_t \gets n_t + 1$.

In each iteration of the main loop we change the number of tied entries, even if 
the iteration resulted in the cone constraint being satisfied and there is no
new tied entry. To compensate for this we let $n_t \gets n_t - 1$ after the
final iteration. (If the cone constraint is satisfied after the first iteration,
then Algorithm \ref{alg:sum-largest} returns $n_t = 1$ even though there are no ties. However, 
the behaviour of the algorithm is still correct in this case since it also returns
$a_t = n_u$.)

What remains is to determine $r$, the ratio of reduction of the untied and tied
entries. Consider fixed values on the untied and tied entries. 
Suppose that we change the untied entries by a quantity $\Delta_u$. The
instantaneuous change of the objective value of the projection problem is
$2 n_u \Delta_u$, and the instantaneuous change in the sum of the $k$ entries 
elements is $n_u$. The ratio of these changes is $2 n_u \Delta_u / n_u =  2
\Delta_u$. Similarily, the instantaneuous change of the objective value of the 
projection problem due to the change in tied entries is $2 n_t \Delta_t$, and
the instantaneuous change in the sum of the $k$ largest entries is $k - n_u$.
The ratio of these changes is $2 n_t \Delta_t/(k - n_u)$. We want to pick $r =
\Delta_u / \Delta_t$ such that the relative changes are equal, \ie 
\[
2 \Delta_u = \frac{2 n_t \Delta_t}{k - n_u}.
\] 
Solving for $r$ gives $r = n_t/(k - n_u)$. Putting together all steps 
above results in Algorithm \ref{alg:sum-largest}.

\newpage
\bibliographystyle{alpha}
\bibliography{references.bib}
\end{document}